\documentclass[aos,MSNbibl,seceqn,nameyear,dvips]{arximspdf}
\usepackage{graphicx}
\doi{10.1214/12-AOS1023} 
\volume{40}
\issue{3}
\pubyear{2012}
\firstpage{1764}
\lastpage{1793}

\makeatletter

\newcommand{\rright}{\right}
\newcommand{\lleft}{\left}
\newcommand{\rrvert}{\vert}
\newcommand{\llvert}{\vert}

\newcommand{\iint}{\int\!\!\int}

\newtheorem{tm}{Theorem}[section]
\newtheorem{pr}[tm]{Proposition}
\newtheorem{lm}[tm]{Lemma}
\newtheorem{kor}[tm]{Corollary}

\newproclaim{rem}[tm]{Remark}

\newcommand{\R}{{\mathbb{R}}}
\newcommand{\N}{{\mathbb{N}}}

\newcommand{\E}{\mathbb{E}}
\renewcommand{\P}{{\mathbb{P}}}

\newcommand{\dif}{\mathrm{d}}

\newcommand{\eps}{\varepsilon}

\newcommand{\one}{\mathbf{1}}

\newcommand{\eqd}{\stackrel{d}{=}}
\newcommand{\dto}{\stackrel{d}{\to}}
\newcommand{\pto}{\stackrel{\P}{\to}}
\newcommand{\argmin}{\mathop{\arg\min}}

\newcommand{\rmd}{\mathrm{d}}
\newcommand{\dw}{\rmd{w}}
\newcommand{\tint}{\int}

\makeatother

\begin{document}
\begin{frontmatter}

\title{An M-estimator for tail dependence in~arbitrary~dimensions}
\runtitle{An M-estimator for tail dependence}

\begin{aug}
\author[A]{\fnms{John H. J.} \snm{Einmahl}\ead[label=e1]{j.h.j.einmahl@uvt.nl}},
\author[B]{\fnms{Andrea} \snm{Krajina}\thanksref{t1}\ead[label=e2]{Andrea.Krajina@mathematik.uni-goettingen.de}}
\and
\author[C]{\fnms{Johan} \snm{Segers}\corref{}\thanksref{t2}\ead[label=e3]{johan.segers@uclouvain.be}}
\runauthor{J. H. J. Einmahl, A. Krajina and J. Segers}
\affiliation{Tilburg University, University of G\"{o}ttingen
and Universit\'e~Catholique~de~Louvain}
\address[A]{J. H. J. Einmahl\\
Department of Econometrics and OR\\
\quad and CentER\\
Tilburg University\\
PO Box 90153\\
5000 LE Tilburg\\
The Netherlands\\
\printead{e1}}
\address[B]{A. Krajina\\
Institute for Mathematical Stochastics\\
University of G\"{o}ttingen, G\"ottingen\\
Germany\\
\printead{e2}}
\address[C]{J. Segers\\
ISBA\\
Universit\'e Catholique de Louvain\\
Voie du Roman Pays, 20\\
B-1348 Louvain-la-Neuve\\
Belgium\\
\printead{e3}} 
\end{aug}

\thankstext{t1}{Supported by the Open Competition grant from the Netherlands
Organisation for Scientific Research (NWO) and from the Deutsches
Forschungsgemeinschaft (DFG) Grant SNF FOR 916.
Research was mainly performed at Tilburg University and Eurandom
(Eindhoven).}

\thankstext{t2}{Supported by IAP research network Grant P6/03
of the Belgian government (Belgian Science Policy) and from the
contract ``Projet d'Actions de Recherche Concert\'ees'' no. 07/12/002
of the Communaut\'e fran\c{c}aise de Belgique,
granted by the Acad\'emie universitaire Louvain.}

\received{\smonth{9} \syear{2011}}
\revised{\smonth{4} \syear{2012}}

%
\begin{abstract}
Consider a random sample in the max-domain of attraction of a
multivariate extreme value distribution such that the dependence
structure of the attractor belongs to a parametric model. A new
estimator for the unknown parameter is defined as the value that
minimizes the distance between a vector of weighted integrals of the
tail dependence function and their empirical counterparts. The
minimization problem has, with probability tending to one, a unique,
global solution. The estimator is consistent and asymptotically normal.
The spectral measures of the tail dependence models to which the method
applies can be discrete or continuous. Examples demonstrate the
applicability and the performance of the method.
\end{abstract}

%

\end{frontmatter}

\section{Introduction}\label{sec1}
\label{S4Intro}
Statistics of multivariate extremes finds important applications in
fields like finance, insurance, environmental sciences, aviation
safety, hydrology
and meteorology. When considering multivariate extreme events, the
estimation of the tail dependence structure is the key part of the statistical
inference. This tail dependence structure, represented by the stable
tail dependence function~$l$, becomes rather complex if the dimension increases.
Therefore, it is customary to model this multivariate function $l$
parametrically, which leads to a semiparametric model.
The interest in parametric tail dependence models has existed since the
early sixties of the 20th century [\citet{Gumbel60}],
but new models are still being proposed [\citet{BD07,CDN10,BS11}]. Most
of the existing estimators of the parameter, $\theta$,
are likelihood-based and their asymptotic behavior is only known in
dimension two [\citet{CT91,JSW92,Smith94,LT96,dHNP08,GPS11}].
For many applications, the bivariate setup is too restrictive. Also,
the likelihood-based estimation methods exclude models that entail a
nondifferentiable function $l$, like the widely used factor models; see
(\ref{E4FM}) below.

It is the goal of this paper to present and provide a comprehensive
treatment of novel M-estimators of $\theta$. The estimators can be used
in arbitrary dimension~$d$.
Moreover, not relying on the differentiability of $l$, the estimators
are broadly applicable. We establish, again for arbitrary dimension
$d$, the asymptotic normality of our estimators, which yields
asymptotic confidence regions and tests for the parameter $\theta$.
The results in this paper make statistical inference possible for many
multivariate extreme value models that either cannot be handled at all
by currently available methods or for which statistical theory has only
been provided for the bivariate case. Monte Carlo simulation studies
confirm that our estimators perform well in practice; see Sections \ref
{S5Logistic} and~\ref{S6FM}.

The present estimators are a major extension of the method of moments
estimators for dimension two [\citet{EKS08}]. For applications, the
crucial difference is that it is now possible to handle truly
multivariate data. Also, theoretically, extreme value analysis in
dimensions larger than two is quite challenging, which explains why in
many papers attention is restricted to the bivariate case. In
particular, we establish the asymptotic behavior of the nonparametric
estimator of $l$ in arbitrary dimensions and under nonrestrictive
smoothness conditions; compare, for instance, with \citet{DH98} in the
bivariate case. Another novel aspect is that the method of moments
technique is replaced by general \mbox{M-estimation}, that is, allowing for
more estimating equations than the dimension of the parameter space.
This more flexible procedure may serve to increase the efficiency of
the estimator.

The absence of smoothness assumptions on $l$ makes it possible to
estimate the tail dependence structure of factor models like
$X=(X_1,\ldots,X_d)$, with
%
%
\begin{equation}
\label{E4FM} X_j=\sum_{i=1}^r
a_{ij}Z_i+\varepsilon_j,\qquad j = 1, \ldots, d,
\end{equation}
consisting of the following ingredients: nonnegative factor loadings
$a_{ij}$ 
and independent, heavy-tailed random variables $Z_i$ called factors;
independent random variables $\varepsilon_j$ whose tails are lighter
than the ones of the factors and which are independent of them. This
kind of factor model is often used in finance, for example, in modeling
market or credit risk [\citet{FF93,MS04,GdHdV07}]. From equation
(\ref
{E4FMl}) below, we see that the stable tail dependence function $l$ of
such a factor model is not everywhere differentiable, causing
likelihood-based methods to break down.

The organization of the paper is as follows. The basics of the
tail dependence structures in multivariate models are presented in
Section~\ref{S2TailDep}. The M-estimator is defined in
Section~\ref{S3Est}. Section~\ref{S4Results} contains the main
theoretical results: consistency and asymptotic normality of the
M-estimator, and some consequences of the asymptotic normality
result that can be used for construction of confidence regions and for
testing. This section also contains the asymptotic normality result for
$\hat{l}_n$.
In Section~\ref{S5Logistic} we apply the M-estimator to
the well-known logistic stable tail dependence function (\ref
{E4logL}). The tail
dependence structure of factor models is studied in
Section~\ref{S6FM}. Both models are illustrated with simulated and
real data. The proofs are deferred to Section~\ref{S7Proofs}.


\section{Tail dependence}\label{sec2}
\label{S2TailDep}

We will write points in $\R^d$ as $x = (x_1,\ldots,x_d)$ and random
vectors as $X_i = (X_{i1},\ldots,X_{id})$, for $i=1,\ldots,n$. Let
$X_1,\ldots, X_n$ be independent random vectors in $\R^d$ with common
continuous distribution function $F$ and marginal distribution
functions $F_1, \ldots, F_d$. %
%
For $j=1,\ldots,d$, write $M_n^{(j)}:= \max_{i=1,\ldots,n}X_{ij}$. We
say that $F$ is in the max-domain of attraction of an extreme value
distribution $G$ if there exist sequences $a_n^{(j)}>0$, $b_n^{(j)}\in
\R$, $j=1,\ldots,d$, such that
%
%
\begin{equation}
\label{eqMDA} \lim_{n\to\infty}\P\biggl( \frac{M_n^{(1)} - b_n^{(1)}
}{a_n^{(1)}}\leq
x_1,\ldots,\frac{M_n^{(d)} - b_n^{(d)}
}{a_n^{(d)}}\leq x_d \biggr) = G(x)
\end{equation}
for all continuity points $x \in\R^d$ of $G$. The margins $G_1,
\ldots, G_d$ of $G$ must be univariate extreme value distributions and
the dependence structure of $G$ is determined by the relation
\[
- \log G(x) = l \bigl( - \log G_1(x_1), \ldots, - \log
G_d(x_d) \bigr)
\]
for all points $x$ such that $G_j(x_j) > 0$ for all $j = 1, \ldots,
d$. The stable tail dependence function $l\dvtx[0, \infty)^d \to[0,
\infty)$ can be retrieved from $F$ via
%
%
\begin{equation}\qquad
\label{E4defl} l(x) = \lim_{t\downarrow0}t^{-1} \P\bigl
\{1-F_1(X_{11})\leq tx_1\mbox{ or }\ldots\mbox{
or }1-F_d(X_{1d})\leq tx_d \bigr\}.
\end{equation}
In fact, the joint convergence in~(\ref{eqMDA}) is equivalent to
convergence of the $d$ marginal distributions together with~(\ref{E4defl}).

In this paper we will only assume the weaker relation~(\ref{E4defl}).
By itself,~(\ref{E4defl}) holds if and only if the random vector $( 1
/ \{ 1 - F_1(X_{1j}) \} )_{j=1}^d$ belongs to the max-domain of
attraction of the extreme value distribution $G_0(x) = \exp\{ -
l(1/x_1, \ldots, 1/x_d) \}$ for $x \in(0, \infty)^d$. Alternatively,
the existence of the limit in~(\ref{E4defl}) is equivalent to
multivariate regular variation of the random vector $( 1 / \{ 1 -
F_1(X_{1j}) \} )_{j=1}^d$ on the cone $[0, \infty]^d \setminus\{(0,
\ldots, 0)\}$ with limit measure or exponent measure $\mu$ given by
%
\[
\mu\bigl( \bigl\{z\in[0,\infty]^d\dvtx z_1
\geq x_1\mbox{ or }\ldots\mbox{ or } z_d\geq
x_d \bigr\} \bigr) = l(1/x_1, \ldots, 1/x_d)
\]
[\citet{R87,BGST04,dHF06}].
The measure $\mu$ is homogeneous, that is, $\mu(tA) = t^{-1}\mu(A)$,
for any $t>0$ and any relatively compact Borel set $A\subset[0,\infty
]^d \setminus\{(0, \ldots, 0)\}$, where \mbox{$tA:= \{tz\dvtx z\in A\}$}. This
homogeneity property yields a decomposition of $\mu$ into a radial and
an angular part [\citet{dHR77,R87}]. Let $\Delta_{d-1}:=\{w \in[0,
1]^d\dvtx w_1+\cdots+w_d=1\}$ be the unit simplex in $\R^d$. Associated
to $B\subset\Delta_{d-1}$ and $t>0$ is the set
\[
B_{t} = \Biggl\{x\in[0,\infty)^d\setminus\bigl\{(0,
\ldots,0)\bigr\} \dvtx{ \sum_{j=1}^d}x_j
\geq t, x \Big/ { \sum_{j=1}^d}x_j
\in B \Biggr\}.
\]
By the homogeneity property of the exponent measure, it holds that $\mu
(B_t) = t^{-1}\mu(B_1)$. Writing $H(B) = \mu(B_1)$ defines a finite
measure $H$ on $\Delta_{d-1}$, called the spectral or angular measure.
Any finite measure satisfying the moment conditions
%
%
\begin{equation}
\label{E4moment} \int_{\Delta_{d-1}}w_jH(\dif w) = 1,\qquad j=1,
\ldots,d,
\end{equation}
is a spectral measure. Adding up the $d$ constraints in
(\ref{E4moment}) shows that $H/d$ is a probability measure.

Sometimes it is more convenient to work with the measure $\Lambda$
obtained from $\mu$ after the
transformation $(x_1,\ldots,x_d)\mapsto(1/x_1,\ldots,1/x_d)$. The
meas\-ure~$\Lambda$ is also called the exponent measure and
it satisfies the homogeneity property $\Lambda(tA) = t\Lambda(A)$,
for any $t>0$ and
Borel set $A\subset[0,\infty]^d\setminus\{(\infty,\ldots,\infty)\}$.


There is a one-to-one correspondence between the stable tail dependence
function $l$, the exponent measures $\mu$ and $\Lambda$, and the
spectral measure $H$. In particular, we have
%
%
\begin{eqnarray}
\label{E4lmu}\qquad
l(x) &=& \mu\bigl( \bigl\{(z_1,\ldots,z_d)\in[0,
\infty]^d\dvtx z_1\geq1/x_1 \mbox{ or }\ldots
\mbox{ or } z_d\geq1/x_d \bigr\} \bigr)
\\
\label{E4lLambda}
&=& \Lambda\bigl( \bigl\{(u_1,\ldots,u_d)\in[0,
\infty]^d \dvtx u_1\leq x_1 \mbox{ or }\ldots
\mbox{ or } u_d\leq x_d \bigr\} \bigr)
\\
\label{E4lspectralm}
&=& \int_{\Delta_{d-1}}\max_{j=1,\ldots,d}\{ w_jx_j
\} H(\dif w).
\end{eqnarray}

From the above representations and the moment constraints (\ref
{E4moment}), it follows that the function $l$ has the following properties:
\begin{itemize}
\item$\max\{x_1,\ldots,x_d\}\leq
l(x)\leq x_1+\cdots+x_d$ for all $x \in[0, \infty)^d$; in particular,
$l(z,0, \ldots, 0) = \cdots= l(0, \ldots, 0,z) = z$ for all $z \geq0$;
\item$l$ is convex; and
\item$l$ is homogeneous of order one: $l(tx_1,\ldots,tx_d)=t
l(x_1,\ldots,x_d)$, for all $t > 0$ and all $x\in[0,\infty)^d$.
\end{itemize}
The function $l$ is connected to the function $V$ in \citet{CT91}
through $l(x) = V(1/x_1, \ldots, 1/x_d)$ for $x \in(0,
\infty)^d$.

The right-hand partial derivatives of $l$ always exist; indeed, by
bounded convergence it follows that for $j=1,\ldots,d$, as $h
\downarrow0$,
%
%
\begin{eqnarray}
\label{E4lgrad}
&&
\frac{1}{h} \bigl( l(x_1,\ldots,
x_{j-1},x_j+h,x_{j+1},\ldots, x_d) -
l(x_1,\ldots, x_{j-1},x_j,x_{j+1},\ldots,
x_d) \bigr)
\nonumber
\\
&&\qquad= \int_{\Delta_{d-1}} \frac{1}{h} \Bigl( \max\Bigl\{
w_j x_j + w_j h, \max_{s\neq j}
\{w_s x_s\} \Bigr\} \nonumber\\[-8pt]\\[-8pt]
&&\hspace*{93.1pt}{} - \max\Bigl\{ w_j
x_j, \max_{s\neq j} \{w_s x_s\} \Bigr
\} \Bigr) H(\dw)
\nonumber\\
&&\qquad\to\int_{\Delta_{d-1}} w_j \mathbf{1} \Bigl\{
w_j x_j \geq\max_{s\neq
j} \{w_s
x_s\} \Bigr\} H(\dw).
\nonumber
\end{eqnarray}
Similarly, the left-hand partial derivatives exist for all $x \in
(0, \infty)^d$. By convexity, the function $l$ is almost everywhere
continuously differentiable, with its gradient vector of
(the right-hand) partial derivatives as in~(\ref{E4lgrad}).

\section{Estimation}\label{sec3}
\label{S3Est}

\label{SEst} Let $R_i^j$ denote the rank of $X_{ij}$ among
$X_{1j},\ldots,X_{nj}$, $i=1,\ldots,n$, $j=1,\ldots,d$. For
$k\in\{1,\ldots,n\}$, define a nonparametric estimator of~$l$ by
%
%
\begin{eqnarray}
\label{E4lnhat} \hat{l}_n(x)&=&\hat{l}_{k,n}(x)\nonumber\\[-8pt]\\[-8pt]
:\!&=&
\frac{1}{k}\sum_{i=1}^n\one\biggl
\{R_{i}^1> n+\frac{1}{2}-kx_1 \mbox{ or }
\ldots\mbox{ or } R_{i}^d>n+\frac{1}{2}-kx_d
\biggr\};\nonumber
\end{eqnarray}
see \citet{Huang92} and \citet{DH98} for the bivariate case. This
definition follows from~(\ref{E4defl}), with all the distribution
functions replaced by their empirical counterparts, and with $t$
replaced by $k/n$. Here $k = k_n$ is such that $k\to\infty$ and
$k/n\to0$ as $n\to\infty$. The constant $1/2$ in the argument of the
indicator function helps to improve the finite-sample properties of the
estimator.

In the literature, the stable tail dependence function
is often modeled parametrically. We impose that the stable tail
dependence function $l$ belongs to some parametric family
$\{l( \cdot; \theta)\dvtx\theta\in\Theta\}$, where $\Theta
\subset\R^p$, $p \geq1$. Note that this is still a large, flexible
model since there is no restriction on the marginal distributions and
the copula is constrained only through $l$; see~(\ref{E4defl}).

We propose an M-estimator of $\theta$.
Let $q\geq p$. Let $g\equiv(g_1,\ldots,g_q)^T\dvtx[0,1]^d\to\R^q$ be a
column vector of integrable functions such that $\varphi\dvtx\Theta\to
\R^q$ defined
by
%
%
\begin{equation}
\label{E4phi} \varphi(\theta):=\int_{[0,1]^d} g(x) l(x;\theta)
\,\dif x
\end{equation}
is a homeomorphism between $\Theta$ and its image $\varphi(\Theta)$.
Let $\theta_0$ denote the true parameter value. The M-estimator
$\hat{\theta}_n$ of $\theta_0$ is defined as a minimizer of the
criterion function
%
%
\begin{equation}\quad
\label{eqQkn} Q_{k,n}(\theta) = \biggl\| \varphi(\theta) - \int g
\hat{l}_n \biggr\|^2 = \sum_{m=1}^q
\biggl( \int_{[0,1]^d}g_m(x) \bigl(
\hat{l}_n(x)-l(x;\theta) \bigr)\,\dif x \biggr)^2,\hspace*{-28pt}
\end{equation}
where \mbox{$\|\cdot\|$} is the Euclidean norm. In other words, if
$\hat{Y}_n={\argmin_{y\in\varphi(\Theta)}}\|y-\int g\hat{l}_n\|$, then
$\hat{\theta}_n\in\varphi^{-1}(\hat{Y}_n)$. Later we show that
$\hat{\theta}_n$ is, with probability tending to one, unique.

The fact that our model assumption only concerns a limit relation in
the tail shows up in the estimation procedure through the choice of
$k$, which determines the effective sample size.
When we study asymptotic properties of either $\hat{l}_n$ or $\hat
{\theta}_n$, $k = k_n$ is an intermediate sequence,
that is, $k \to\infty$ and $k/n \to0$ as $n\to\infty$. In
practice, the choice of optimal $k$ is a notorious problem, and here we
address this issue in the usual way: we
present the finite sample results over a wide range of $k$; see
Sections~\ref{S5Logistic} and~\ref{S6FM}.

%
\begin{rem}
The estimator $\hat\theta_n$ depends on $g$. In line with the
classical method of moments and for computational feasibility, we will
choose $g$
to be a vector of low degree polynomials.
In Sections~\ref{sec5} and~\ref{sec6} we will see that the obtained
estimators have a
good performance and a wide applicability. Finding an
optimal $g$ is very difficult and statistically not very
useful since such a $g$ depends on the true, unknown $\theta_0$. For
example, when \mbox{$p=q=1$}, a~function $g$ that minimizes the asymptotic
variance is $(\partial/\partial\theta) l(x;\theta_0)$.
For two-dimensional and five-dimensional data, a sensitivity analysis
on the choice of~$g$ is performed in Section~\ref{sec5}. Simple
functions like
$1$ or $x_1$ lead to
estimators that perform approximately the same as the pseudo-estimator
based on the optimal $g$. This supports our choices of $g$ and also
suggests that the estimator is not so sensitive to the
choice of $g$.
\end{rem}

%
\begin{rem}
Since $l$, part of the model, is parametrically specified, in
principle, pseudo maximum likelihood estimation could be used. This
method, however, does not apply to many
interesting models where $l$ is not differentiable, like the factor
model in~(\ref{E4FM}).
Moreover, no theory is known for dimensions higher than 2, unless the
limit relation~(\ref{E4defl}) is replaced by an equality for all
sufficiently small~$t$. In this paper, the emphasis is on higher
dimensions and for a large part on the factor model. Therefore, the
pseudo MLE is not an available competitor.\vadjust{\goodbreak}
\end{rem}

\section{Asymptotic results}\label{sec4} \label{SResults}
\label{S4Results}
Let $\hat{\Theta}_n$ be the set of minimizers of $Q_{k,n}$ in (\ref
{eqQkn}), that is,
\[
\hat{\Theta}_n:= \argmin_{\theta\in\Theta} \biggl\| \varphi(\theta) - \tint g
\hat{l}_n \biggr\|^2.
\]
Note that $\hat{\Theta}_n$ may be empty or may contain more than one
element. We show that under suitable conditions, a minimizer exists,
that it is unique with probability tending to one, and that it is a
consistent and asymptotically normal estimator of $\theta_0$. In addition,
we show that the nonparametric estimator $\hat{l}_n$ in (\ref
{E4lnhat}) is asymptotically normal.

\subsection{Notation}\label{sec4.1}

Recall the definition of the measure $\Lambda$ from Section \ref
{S2TailDep}. Let $W_\Lambda$ be a mean-zero Wiener process indexed by
Borel sets
of $[0,\infty]^d\setminus\{(\infty,\ldots,\infty)\}$ with ``time''
$\Lambda$: its covariance structure is given by
%
%
\begin{equation}
\label{E4ELambda}\E\bigl[ W_\Lambda(A_1) W_\Lambda(A_2)
\bigr] = \Lambda(A_1 \cap A_2)
\end{equation}
for any two Borel sets $A_1$ and $A_2$ in $[0,\infty]^d\setminus\{
(\infty,\ldots,\infty)\}$. Define
%
%
\begin{equation}
\label{defW} W_l(x):= W_\Lambda\bigl(\bigl\{u\in[0,
\infty]^d\setminus\bigl\{(\infty,\ldots,\infty)\bigr\}\dvtx u_1
\leq x_1 \mbox{ or }\ldots\mbox{ or } u_d\leq
x_d\bigr\}\bigr).\hspace*{-35pt}
\end{equation}
Let $W_{l,j}, j=1,\ldots,d$, be the marginal processes
%
%
\begin{equation}
\label{defWj} W_{l,j}(x_j):= W_l(0,
\ldots,0,x_j,0,\ldots,0), \qquad x_j\geq0.
\end{equation}
Define $l_j$ to be the right-hand partial derivative of $l$ with
respect to $x_j$, where $j=1,\ldots,d$ [see~(\ref{E4lgrad})]; if $l$
is differentiable, $l_j$ is equal to the corresponding partial
derivative of $l$. Write
%
%
\begin{equation}\qquad
\label{eqB} B(x):= W_l(x)-\sum_{j=1}^d
l_j(x) W_{l,j}(x_j), \qquad\tilde{B}:=\int
_{[0,1]^d} g(x) B(x) \,\dif x.
\end{equation}
The distribution of $\tilde{B}$ is zero-mean Gaussian with
covariance matrix
%
%
\begin{equation}
\label{E4Sigma} \Sigma:=\iint_{([0,1]^d)^2} \E\bigl[B(x) B(y)\bigr] g(x)
g(y)^T \,\dif x \,\dif y \in\R^{q \times q}.
\end{equation}
Note that if $l$ is parametric, $\Sigma$ depends on the parameter,
that is, $\Sigma= \Sigma(\theta)$.

Assuming $\theta$ is an interior point of $\Theta$ and $\varphi$ is
differentiable in $\theta$, let $\dot{\varphi}(\theta) \in\R^{q\times
p}$ be the
total derivative of $\varphi$ at $\theta$, and, provided $\dot
{\varphi}(\theta)$ is of full rank, put
%
%
\begin{equation}
\label{E4matrixM} M(\theta):= \bigl(\dot{\varphi}(\theta)^T\dot{
\varphi}(\theta) \bigr)^{-1}\dot{\varphi}(\theta)^T \Sigma(
\theta) \dot{\varphi}(\theta) \bigl(\dot{\varphi}(\theta)^T\dot{
\varphi}(\theta) \bigr)^{-1}\in\R^{p\times p}.
\end{equation}

\subsection{Results}\label{sec4.2}

We state the asymptotic results for the M-estimator, $\hat{\theta
}_n$, and the asymptotic normality of $\hat{l}_n$. The latter is a
result of independent interest, and requires continuous
partial derivatives of $l$, which is not an assumption for the
asymptotic normality of the M-estimator. The proofs can be found in
Section~\ref{S7Proofs}.\vadjust{\goodbreak}

%
\begin{tm}[(Existence, uniqueness and consistency of $\hat{\theta
}_n$)]\label{Tcons}
Let $g\dvtx\break[0,1]^d \to\R^q$ be integrable.
\begin{longlist}[(ii)]
\item[(i)] If $\varphi$
is a homeomorphism from $\Theta$ to $\varphi(\Theta)$ and if there exists
$\varepsilon_0>0$ such that the set $\{\theta\in\Theta\dvtx
\|\theta-\theta_0\|\leq\varepsilon_0\}$ is closed, then for
every $\eps$ such that $\eps_0\geq\eps> 0$, as $n\to\infty$,
\[
\P\bigl( \hat{\Theta}_n \neq\varnothing\mbox{ and } \hat{
\Theta}_n\subset\bigl\{\theta\in\Theta\dvtx\|\theta-\theta_0\|
\leq\varepsilon\bigr\} \bigr) \to1.
\]
\item[(ii)] If in addition to the assumptions of \textup{(i)}, $\theta_0$ is
in the interior of the parameter space, $\varphi$ is twice
continuously differentiable and $\dot{\varphi}(\theta_0)$ is of
full rank, then, with probability tending to one, $Q_{k,n}$ in (\ref
{eqQkn}) has a unique minimizer $\hat{\theta}_n$. Hence,
\[
\hat{\theta}_n\pto\theta_0\qquad\mbox{as } n\to\infty.
\]
\end{longlist}
\end{tm}
In part (i) of this theorem we assume that the set $\{\theta\in\Theta
\dvtx\|\theta-\theta_0\|\leq\varepsilon\}$ is closed for some
$\varepsilon>0$. This is a generalization of the usual assumption that
$\Theta$ is open or closed, and includes a wider range of possible
parameter spaces.

%
\begin{tm}[(Asymptotic normality of $\hat{\theta}_n$)]\label{Tantheta}
If in addition to the assumptions of Theorem~\ref{Tcons}\textup{(ii)},
the following two conditions hold:
\begin{longlist}[(C2)]
\item[(C1)] $t^{-1}\P\{1-F_1(X_{11})\leq tx_1\mbox{ or }\ldots\mbox
{ or }1-F_d(X_{1d})\leq
tx_d\}-l(x)=O(t^\alpha)$, uniformly in $x\in\Delta_{d-1}$ as
$t\downarrow0$,
for some $\alpha>0$,\vspace*{1pt}
\item[(C2)] $k=o (n^{2\alpha/(1+2\alpha)} )$, for the
positive number $\alpha$ of \textup{(C1)}, and $k\to\infty$ as $n\to\infty$,
\end{longlist}
then as $n\to\infty$, with $M$ as in~(\ref{E4matrixM}),
%
%
\begin{equation}
\label{E4ANtheta} \sqrt{k}(\hat{\theta}_n-\theta_0)\dto
N \bigl( 0,M(\theta_0) \bigr).
\end{equation}
\end{tm}

The following consequence of Theorem~\ref{Tantheta} can be used
for the construction of confidence regions. Recall from (\ref
{E4lspectralm}) that $H_\theta$ is the spectral measure corresponding
to $l(\cdot;\theta)$. Let $\chi^2_\nu$ denote the $\chi^2$-distribution
with $\nu$ degrees of freedom.

%
\begin{kor}\label{Tcorollary}
If in addition to the conditions of Theorem~\ref{Tantheta}, the
map $\theta\mapsto H_\theta$ is weakly continuous at $\theta_0$
and if the matrix $M(\theta_0)$ is nonsingular, then
as $n\to\infty$,
%
%
\begin{equation}
\label{E4confreg} k(\hat{\theta}_n-\theta_0)^T
M(\hat{\theta}_n)^{-1}(\hat{\theta}_n-
\theta_0) \dto\chi^2_p.
\end{equation}
\end{kor}

Let $1\leq r<p$ and $\theta= (\theta_1,
\theta_2)\in\Theta\subset\R^p$, where $\theta_1\in\R^{p-r}$,
$\theta_2\in\R^r$. We want to test $\theta_2=\theta_{2}^*$ against
$\theta_2\neq\theta_{2}^*$, where $\theta_2^*$ corresponds to a
submodel. Denote $\hat{\theta}_n=(\hat{\theta}_{1n},
\hat{\theta}_{2n})$, and let $M_2(\theta)$ be the $r\times r$ matrix
corresponding\vadjust{\goodbreak} to the lower right corner of $M$, as below:
%
%
\begin{equation}
\label{E4MandM2} M = \pmatrix{
\cdots\hspace*{2pt}|\hspace*{2pt} \cdots
\cr
\hline
\cdots\hspace*{2pt}|\hspace*{1pt} M_2
}
\in
\R^{p\times p}.
\end{equation}

%
\begin{kor}[(Test)]\label{TTest}
If the assumptions of Corollary~\ref{Tcorollary} are satisfied, and
$\theta_0=(\theta_1,\theta_2^*)\in\Theta$ for some $\theta_1$, then
as $n\to\infty$,
%
%
\begin{equation}
\label{E4test} k\bigl(\hat{\theta}_{2n}-\theta_2^*
\bigr)^TM_2\bigl(\hat{\theta}_{1n},
\theta_2^*\bigr)^{-1}\bigl(\hat{\theta}_{2n}-
\theta_2^*\bigr)\dto\chi_r^2.
\end{equation}
\end{kor}

The above result can be used for testing for a submodel. For
example, we could test for the symmetric logistic model of (\ref
{EsymmLog}) within the
asymmetric logistic one; see Section~\ref{S5Logistic}.

%
\begin{rem}
The matrices $M$ and $M_2$ are needed for the computation of the
confidence regions and the test
statistics. However, computing these matrices can be challenging. To
compute $M$, we first need the $q \times p$ matrix~$\dot{\varphi
}(\theta)$, whose $(i,j)$th element is given by
$\tint g_i(x)(\partial/\partial\theta_j)l(x;\theta)\,\dif x$. The
expression itself will depend on the model in use, but usually the
(right-hand) partial derivatives of $l$ can be computed explicitly,
whereas the integral is to be computed numerically in most cases.
Second, we need to calculate the covariance of the process $\tilde
{B}$. We see from~(\ref{E4Sigma}) that the most difficult part will be
the expression $\E[B(x)B(y)]$. It holds that
\begin{eqnarray*}
\E\bigl[B(x)B(y)\bigr] &=& \E\bigl[W_l(x)W_l(y)\bigr] -
\sum_{j=1}^d l_j(y)\E
\bigl[W_l(x) W_{l,j}(y_j)\bigr]\\[-3pt]
&&{}- \sum
_{i=1}^d l_i(x)\E\bigl[W_{(l,i)}(x_i)
W_l(y)\bigr]
\\[-3pt]
&&{}+\sum_{i=1}^d \sum
_{j=1}^d l_i(x) l_j(y)\E
\bigl[W_{(l,i)}(x_i) W_{l,j}(y_j)\bigr].
\end{eqnarray*}
Using~(\ref{E4ELambda}),~(\ref{defW}),~(\ref{defWj}) and the
relation between $\Lambda$ and $l$, we can express this in $l$ and its
partial derivatives. Numerical integration is then performed
to obtain $\Sigma$.\vspace*{-3pt}
\end{rem}

Finally, we show the asymptotic normality of $\hat{l}_n$. This result
is of independent interest and can be found in the literature for $d=2$
only and under stronger smoothness conditions
on $l$; see \citet{Huang92,DH98} and \citet
{dHF06}. Here, a
large part of its proof is necessary for the proof of the asymptotic
normality of $\hat{\theta}_n$, but we wish to emphasize that the
asymptotic normality of $\hat{\theta}_n$ holds \textit{without} any
differentiability conditions on $l$. Note that under assumption (C3) below,
the process~$B$ in~(\ref{eqB}) is continuous, although $l_j$ may be
discontinuous at points $x$ such that $x_j = 0$.\vadjust{\goodbreak}

The result is stated in an approximation setting, where $\hat{l}_n$
and $B$ are defined on the same probability space obtained by a
Skorohod construction. The random quantities involved are only in
distribution equal to the original ones, but for convenience this is
not expressed in the notation.

%
\begin{tm}[(Asymptotic normality of $\hat{l}_n$ in arbitrary
dimensions)]\label{TanL} If in addition to the conditions \textup{(C1)}
and \textup{(C2)} from Theorem \ref {Tantheta}, the following condition
holds:

\begin{longlist}[(C3)]
\item[(C3)] for all $j=1, \ldots, d$, the first-order partial
derivative of $l$ with respect to $x_j$ exists and is continuous on the
set of points $x$ such that $x_j > 0$,
\end{longlist}
then for every $T>0$, as $n\to\infty$,
%
%
\begin{equation}
\label{E4ANl} \sup_{x\in[0,T]^d}\bigl\llvert\sqrt{k} \bigl(
\hat{l}_n(x)-l(x) \bigr)-B(x)\bigr\rrvert\pto0.
\end{equation}
\end{tm}

\section{Example 1: Logistic model}\label{sec5}
\label{S5Logistic}

The multivariate logistic distribution function with standard Fr\'echet
margins is defined by
\[
F(x_1,\ldots,x_d; \theta)=\exp\Biggl\{- \Biggl({ \sum
_{j=1}^d} x_j^{-1/\theta}
\Biggr)^\theta\Biggr\}
\]
for $x_1>0,\ldots,x_d>0$ and $\theta\in[0,1]$, with the proper limit
interpretation for $\theta= 0$. The corresponding stable tail
dependence function is given by
%
%
\begin{equation}
\label{E4logL} l(x_1,\ldots,x_d; \theta)=
\bigl(x_1^{1/\theta}+\cdots+x_d^{1/\theta
}
\bigr)^\theta.
\end{equation}
Introduced in \citet{Gumbel60}, it is one of the oldest parametric
models of tail dependence.

\subsection*{Sensitivity analysis}
Here we observe how for the logistic model the M-estimator changes with
different choices of $k$, and for different functions $g$. Within this
model, $p=1$ and in the simple case of $p=q=1$, it
is easy to see that the optimal choice for the function $g$ is
$(\partial/\partial\theta)l(x;\theta_0)$. Since it depends on the
unknown true parameter, this is not a viable option for use in practice,
but, as demonstrated below, some simple alternatives result in
estimators with basically the same finite-sample behavior.

%
\begin{figure}

\includegraphics{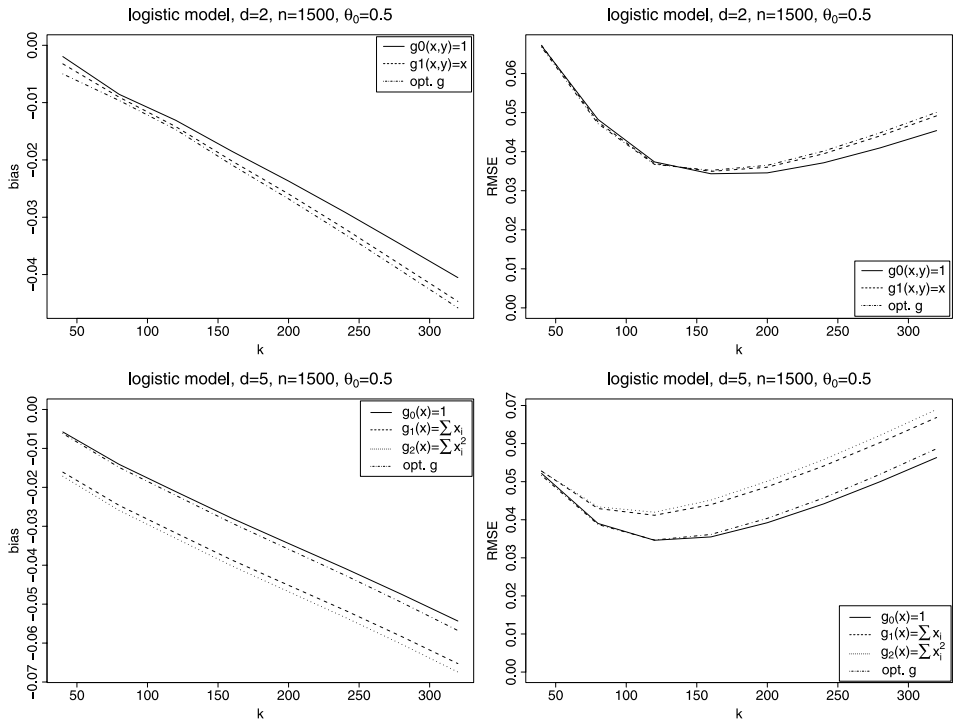}

\caption{Logistic model: the M-estimator for different functions $g$
in dimension $d=2$ (top) and $d=5$ (bottom).}
\label{FigLogD2-5sens}
\end{figure}

The following analysis is performed for the logistic model with $\theta
_0=0.5$, in dimensions 2 and 5. For both settings, we look at 200
replications of samples of size $n=1500$, and take the threshold
parameters $k\in\{40,80,\ldots,320\}$. In the bivariate case we
compare $g_0(x_1,x_2)=1$, $g_1(x_1,x_2)=x_1$ and $g_{\mathrm
{opt}}(x_1,\break x_2)= (\partial/\partial\theta)l(x_1,x_2;\theta_0)$ as
choices for $g$. In the five-dimensional case the functions $g_0$ and
$g_{\mathrm{opt}}$ are defined analogously, and we compare them to
two other functions, $g_1(x)=\sum_{j=1}^5 x_j$ and $g_2(x)=\sum_{j=1}^5
x_j^2$. We use\vspace*{1pt} the bias and the Root Mean Squared Error (RMSE)
to assess the performance of the estimators.\vadjust{\goodbreak} The results are presented
in Figure~\ref{FigLogD2-5sens} for dimensions $d=2$ (top) and $d=5$
(bottom). All of the above choices for $g$ result in similar
finite-sample behavior of the estimator, but the simpler function $g$
leads to a somewhat better performance. The RMSEs for some of these $g$
are even lower than the one for $g_{\mathrm{opt}}$, since they yield a
smaller bias.

Based on these findings, for the logistic model in dimensions $2$ and
$5$, we advise the use of the simplest choice of $g$ given by
$g_0(x)=1$, for all $x\geq0$. The choice of $k$ is slightly more
delicate, but it seems that for $n=1500$ in dimensions 2 and~5, the
choices of $k=150$ and $k=100$, respectively, are reasonable.\vspace*{-3pt}

\subsection*{Comparison with maximum likelihood based estimators}
For $d = 2$, we also compare the M-estimator with $g\equiv1$ with the
censored maximum likelihood method [see \citet{LT96}] and with the maximum
likelihood estimator introduced in \citet{dHNP08}. The latter two we
will call the censored MLE and the dHNP MLE, respectively.
For 200 samples, we compute the censored MLE using the function
\texttt{fitbvgpd} from the R package \texttt{POT} [see \citet
{POT}]; the dHNP
MLE is calculated as described in the original article. Since the
thresholds used in these two methods differ,\vadjust{\goodbreak} and since for a different choice
of threshold we get a different estimator, the comparison is not
straightforward. We consider\ the M-estimator and the dHNP MLE over the range
of $k$ values as used above, and for the censored MLE we take the
thresholds such that the expected number of joint exceedances is between
%
%
\begin{figure}

\includegraphics{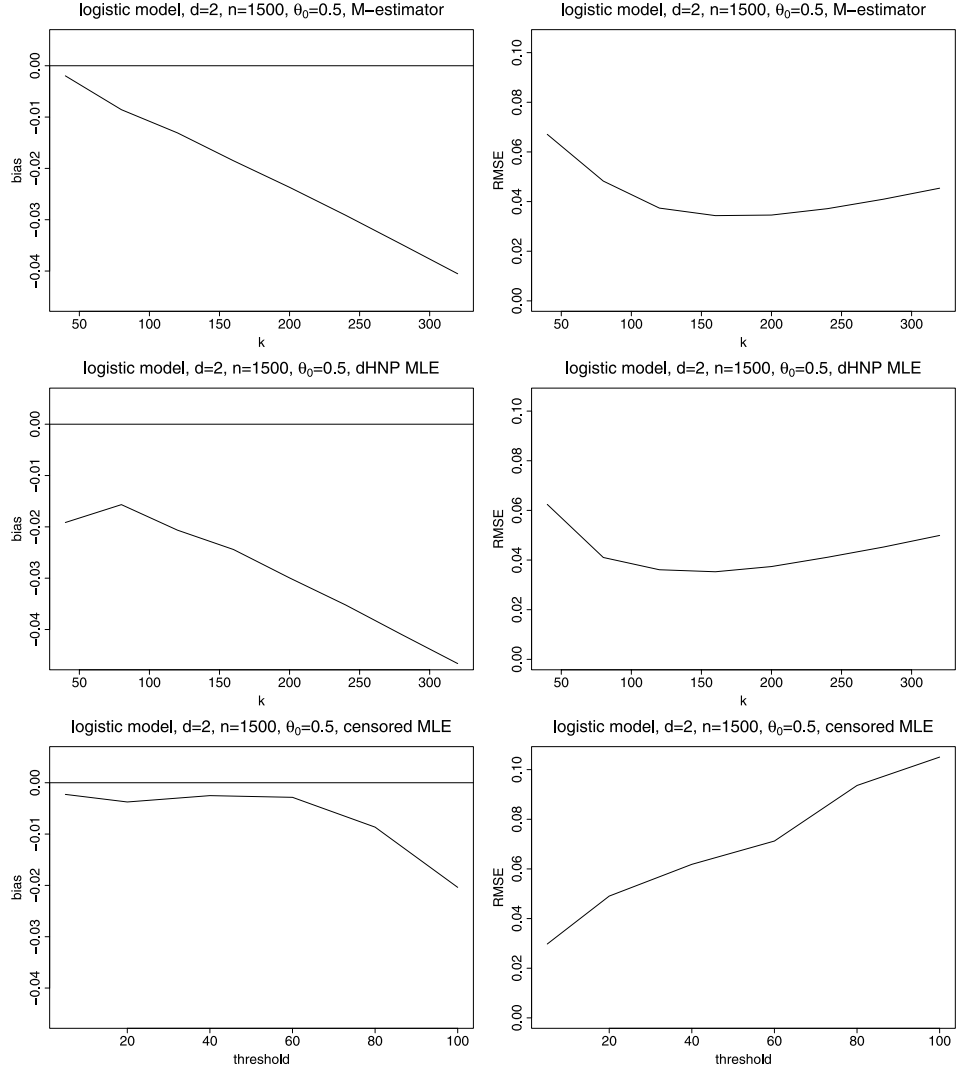}

\caption{The M-estimator with $g(x,y)=g_0(x,y)=1$, the MLE from
de~Haan, Neves and Peng (\protect\citeyear{dHNP08}) and the censored MLE, $d=2$.}
\label{FigLogD2sensML}\vspace*{-3pt}
\end{figure}
$10$ and $160$, approximately, which amounts to thresholds between $5$
and $100$. This way we observe all estimators for their best region
of thresholds. In Figure~\ref{FigLogD2sensML} we see that the methods
perform roughly the same, the RMSEs being of the same order. The lowest
RMSE of the censored MLE ($0.030$) is slightly smaller than the lowest
RMSE of the M-estimator ($0.034$) and
the lowest RMSE of the dHNP estimator ($0.035$), but the M- and the
dHNP estimators are much more robust to the choice of the
threshold.

\subsection*{Further simulation results} 
\begin{figure}
\begin{tabular}{@{}c@{\hspace*{5pt}}c@{}}

\includegraphics{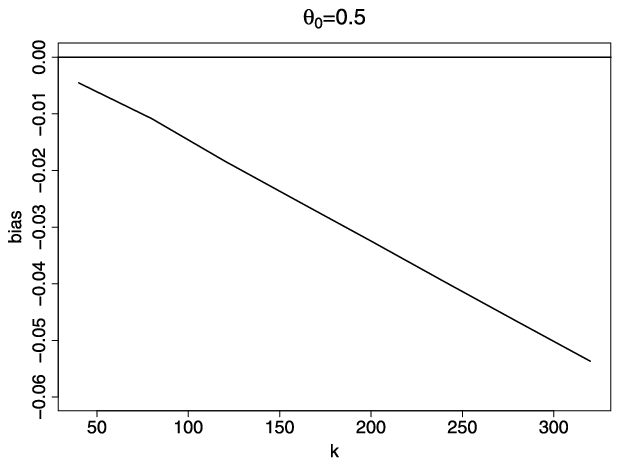}
 & \includegraphics{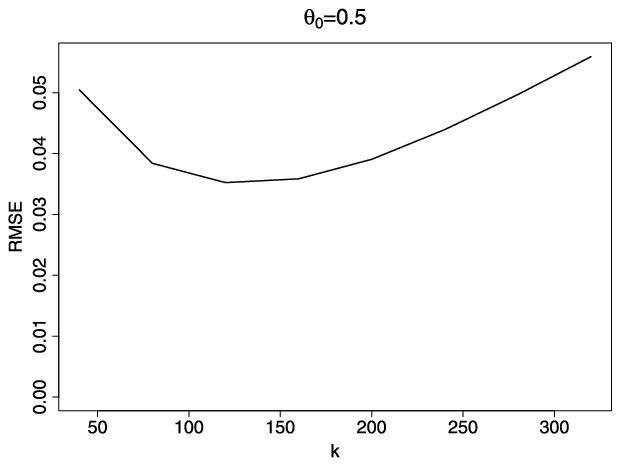}\\
(a) & (b)\\[6pt]

\includegraphics{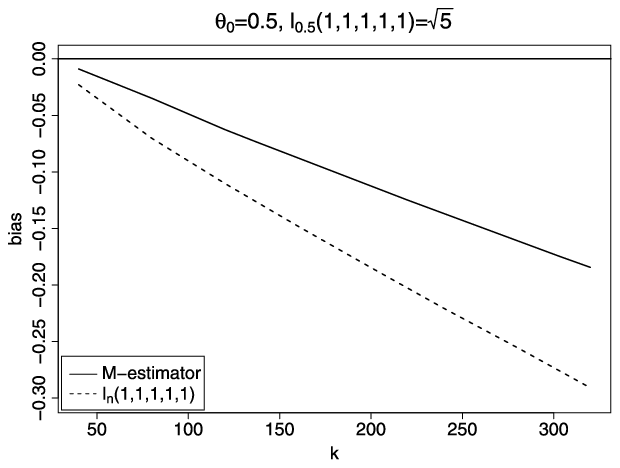}
 & \includegraphics{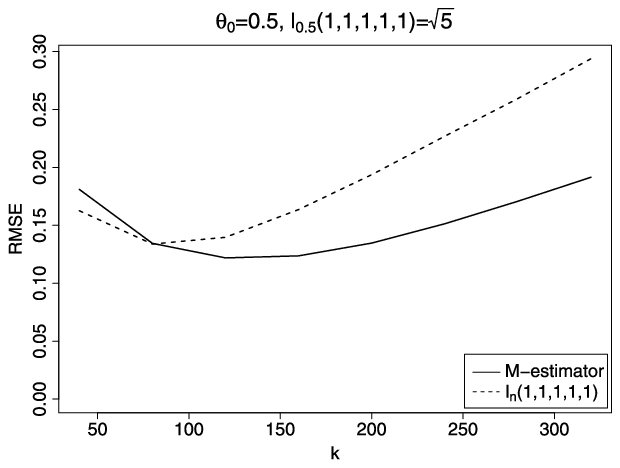}\\
(c) & (d)
\end{tabular}
\caption{Logistic model, $d=5$, $\theta_0=0.5$, $l(1,1,1,1,1; \theta
_0)=\sqrt{5}$. \textup{(a)} Bias of the M-estimator of $\theta$;
\textup{(b)} RMSE of the M-estimator of $\theta$;
\textup{(c)} bias of the estimators of $l(1, 1, 1, 1, 1; 0.5)$;
\textup{(d)} RMSE of the estimators of $l(1, 1, 1, 1, 1; 0.5)$.}
\label{Figlog1}
\end{figure}

We simulate $500$ samples of
size $n=1500$ from a five-dimensional
logistic distribution function with $\theta_0=0.5$. As suggested by
the sensitivity analysis, we opt for $g\equiv1$ when defining $\hat
{\theta}_n$. The bias and the RMSE of this estimator are shown
in the upper panels of Figure~\ref{Figlog1}.

%

Also, we consider the estimation of $l(1,1,1,1,1;\theta)$, based
on this M-estima\-tor~$\hat{\theta}_n$. From~(\ref{E4logL}) it follows
that $l(1,1,1,1,1;\theta) = 5^{\theta}$. The estimator of this
quantity is then $5^{\hat{\theta}_n}$. Since
$\theta_0=0.5$, the true parameter is $\sqrt{5}$.
We compare the bias and the RMSE of this estimator and of the
nonparametric estimator $\hat{l}_n(1,1,1,1,1)$; see~(\ref{E4lnhat}).
The lower panels in Figure~\ref{Figlog1} show that the M-estimator
performs better than
the nonparametric estimator for almost every $k$.\vadjust{\goodbreak}

\subsection*{Real data: Testing and estimation}
\begin{figure}

\includegraphics{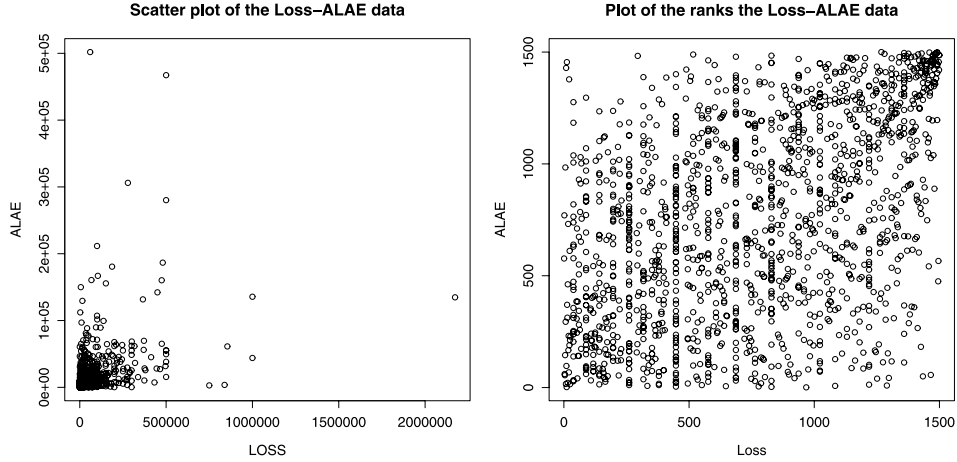}

\caption{The insurance claims Loss-ALAE data.}
\label{FigLossALAE}
\end{figure}

We use the bivariate Loss-ALAE data set, consisting of 1500 insurance
claims, comprising losses and allocated loss adjustment expenses; for
more information, see \citet{FV98}. The scatterplots of
the data and their joint ranks are shown in Figure~\ref{FigLossALAE}.
We consider the asymmetric logistic model described below for their
tail dependence function and
we test whether a more restrictive, symmetric logistic model suffices
to describe the tail dependence of these data.
%
%
The asymmetric logistic tail dependence function was introduced in
\citet{Tawn88} as an extension of the logistic model. In dimension
$d=2$ it is given by
%
%
\begin{equation}
\label{EasymmLog}\qquad l(x,y;\theta, \psi_1,\psi_2) = (1-
\psi_1)x + (1-\psi_2)y+ \bigl((\psi_1x)^{1/\theta}+(
\psi_2y)^{1/\theta} \bigr)^\theta
\end{equation}
with the dependence parameter $\theta\in[0,1]$ and the asymmetry
parameters $\psi_1, \psi_2\in[0,1]$. This model yields a spectral
measure $H$ with atoms at $(1,0)$ and $(0,1)$ whenever $\psi_1 < 1$
and $\psi_2 < 1$. When $\psi_1=\psi_2=:\psi$, we have the symmetric
tail dependence function
%
%
\begin{equation}
\label{EsymmLog} l(x,y;\theta, \psi) = (1-\psi) (x + y) + \psi
\bigl(x^{1/\theta
}+y^{1/\theta} \bigr)^\theta.
\end{equation}

For the given data, we test whether the use of this symmetric model
is justified, as opposed to the wider asymmetric logistic model.
Setting $\eta_1:= (\psi_1+\psi_2)/2 \in[0, 1]$ and $\eta_2:=
(\psi_1-\psi_2)/2 \in[-1/2, 1/2]$, we reparametrize the model in
(\ref{EasymmLog}) so that testing for symmetry amounts to testing
whether $\eta_2=0$.
By Corollary~\ref{TTest}, the test statistic is given by
\[
S_n:= \frac{k \hat{\eta}_{2}^2}{M_2(\hat{\theta},\hat{\eta}_1,0)}.
\]
Table~\ref{TabM-AlogLossALAE1} below shows the obtained values of $S_n$ for the Loss-ALAE
data for selected values of $k$.

\begin{table}
\tablewidth=310pt
\caption{Values of the test statistic $S_n$ for the Loss-ALAE
data for selected values of $k$}
\label{TabM-AlogLossALAE1}
\begin{tabular*}{\tablewidth}{@{\extracolsep{\fill}}lccccc@{}}
\hline
$k$&$50$&$100$&$150$&$200$&$250$\\
$S_n$&$0.041$&$ 0.139 $&$0.294$&$ 0.477 $&$0.681$\\
\hline
\end{tabular*}
\end{table}

Since the critical value is $3.84$, the null hypothesis is clearly not
rejected. Hence, we adopt the symmetric tail dependence model (\ref
{EsymmLog}) and
we compute the M-estimates of $(\theta,\eta_1)=(\theta,\psi)$, the
auxiliary functions
being $g_1(x,y)=x$ and $g_2(x,y)=2(x+y)$. For $k=150$, we obtain $(\hat
{\theta},\hat{\psi})=(0.65, 0.95)$
with estimated standard errors $0.032$ for $\hat{\theta}$ and $0.014$
for $\hat{\psi}$.

\section{Example 2: Factor model}\label{sec6}
\label{S6FM}

Consider the $r$-factor model, $r\in\N$, in dimension $d$:
$X' = (X_1', \ldots, X_d')$ and
%
%
\begin{equation}
\label{EsumFM} X_j' = \sum
_{i=1}^r a_{ij} Z_i+
\varepsilon_j,\qquad j \in\{1, \ldots, d\},
\end{equation}
with $Z_i$ independent Fr\'echet($\nu$) random variables, $\nu> 0$,
with $\varepsilon_j$ independent random variables which have a lighter
right tail than the factors and are independent of them, and with
$a_{ij}$ nonnegative constants such that $\sum_j a_{ij} > 0$ for all
$i$. Factor models of this type are common in various applications;
for example, in finance, see \citet{FF93,MS04,GdHdV07}. However, for
the purpose of studying the tail properties, it is more convenient to consider
the (max) factor model: $X = (X_1, \ldots, X_d)$ and
%
%
\begin{equation}
\label{EmaxFM} X_j = \max_{i=1,\ldots,r} \{a_{ij}
Z_i\},\qquad j \in\{1, \ldots, d\},
\end{equation}
with $a_ {ij}$ and $Z_i$ as above. Note that $X'$ and $X$ have the same
tail dependence function $l$; this essentially follows from the fact
that the ratio of the probabilities of the sum and the maximum of the
$a_{ij}Z_i$ exceeding $x$ tends to 1 as $x \to\infty$ [\citet{EKM97},
page 38].
Let $W_i=Z_i^\nu$, $i=1,\ldots,r$, and
observe that the $W_i$
are standard Fr\'echet random variables. Define a $d$-dimensional
random vector $Y=(Y_1, \ldots,Y_d)$ by
\[
Y_j:= X_j^\nu= \max_{i=1,\ldots,r} \bigl
\{a_{ij}^\nu W_i\bigr\},\qquad j \in\{ 1, \ldots, d\}.
\]
It is easily seen that, as $x\to\infty$,
\[
1 - F_{Y_j}(x) = 1 - \exp\biggl\{-\frac{\sum_{i=1}^r a_{ij}^\nu
}{x} \biggr\}\sim
\frac{\sum_{i=1}^r a_{ij}^\nu}{x}.
\]

Since the $X_j$ variables are increasing transformations of the $Y_j$
variables, the (tail)
dependence structures of $X$ and $Y$ coincide. We will determine
the tail dependence function $l$ and the spectral measure $H$ of
$X$.
%
%
\begin{lm} \label{TlFMl}
Let $X$ follow a factor model given by~(\ref{EsumFM}) or (\ref
{EmaxFM}). Then its stable tail dependence function is given by
%
%
\begin{equation}
\label{E4FMl} l(x_1, \ldots, x_d) = \sum
_{i=1}^r \max_{j=1,\ldots,d} \{b_{ij}
x_j\},\qquad (x_1, \ldots, x_d) \in[0,
\infty)^d,
\end{equation}
where $b_{ij}:=a_{ij}^\nu/\sum_{i=1}^r
a_{ij}^\nu$.
\end{lm}

Next, we are looking for a measure $H$ on the unit simplex
$\Delta_{d-1} = \{ w \in[0, \infty)^d\dvtx w_1+\cdots+w_d = 1 \}$ such
that for all $x \in[0, \infty)^d$,
\[
\sum_{i=1}^r \max_{j=1,\ldots,d}
\{b_{ij} x_j\} = l(x_1, \ldots, x_d)
= \int_{\Delta_{d-1}} \max_{j=1,\ldots,d} \{w_j
x_j\} H(\dif w).
\]
This $H$ is a discrete measure with $r$ atoms given by
%
%
\begin{equation}
\label{atoms} \biggl( \frac{b_{i1}}{\sum_j b_{ij}}, \ldots, \frac
{b_{id}}{\sum_j b_{ij}} \biggr),\qquad i
\in\{1, \ldots, r\},
\end{equation}
the atom receiving mass $\sum_j b_{ij}$, which is positive by
assumption. Such measure $H$ is indeed a spectral measure, for
%
%
\begin{equation}
\label{E4momentH} \int_{\Delta_{d-1}} w_j H(\dif w) = \sum
_{i=1}^r b_{ij} = 1,\qquad j \in\{1,
\ldots, d\}.
\end{equation}
Every discrete spectral measure can arise in this way. This model for
tail dependence is considered also in \citet{LT98}. Extensions to
random fields are considered, for instance, in \citet{WS11}.

The spectral measure is completely determined by the $r\times d$ parameters~$b_{ij}$,
but by the $d$ moment conditions from~(\ref{E4momentH}),
the actual number of parameters is $p = (r-1) d$. The parameter vector
$\theta\in\R^p$, which is to be estimated, can be constructed in
many ways. For identification purposes, the definition of $\theta$
should be unambiguous. We opt for the following approach. Consider
the matrix of the coefficients $b_{ij}$,
\[
\lleft( %
\begin{array} {ccc} b_{11} & \cdots&
b_{r1}
\\
\vdots& \ddots& \vdots
\\
b_{1d} & \cdots& b_{rd}
\end{array}
\rright)\in
\R^{d\times r}.
\]
The coefficients corresponding to the $i$th factor, $i=1,\ldots,r$, are
in the $i$th column of this matrix. We define $\theta$
by stacking the above columns in decreasing order of their sums,
leaving out the column with the lowest sum. (If two columns have the
same sum, we order them then in decreasing order
lexicographically.)

The definition of the M-estimator of $\theta$ involves integrals of
the form
\[
\int_{[0, 1]^d} g_m(x) l(x) \,\dif x = \sum
_{i=1}^r \int_{[0, 1]^d}
g_m(x) \max_{j=1,\ldots,d} \{b_{ij} x_j\} \,\dif
x,\vadjust{\goodbreak}
\]
where $g_m\dvtx[0, 1]^d \to\mathbb{R}$ is integrable and
$m=1,\ldots,q$. A possible choice is $g_m(x) = x_k^s$, where $k \in\{
1, \ldots,
d\}$ and $s \geq0$.
%
%
\begin{lm}\label{TlFMintegrals}
If $l$ is the tail dependence function of a factor model such that all
$b_{ij}>0$, then
%
%
\begin{eqnarray}
\label{Elm62}\qquad
&&
\int_{[0, 1]^d} x_k^s l(x)
\,\dif x\nonumber\\[-8pt]\\[-8pt]
&&\qquad = \sum_{i=1}^r \sum
_{j=1}^d \frac{b_{ij}}{1 + s (1-\delta_{jk})} \int
_0^1 \biggl( \frac{b_{ij}}{b_{ik}} x \wedge1
\biggr)^s \prod_{l=1}^d
\biggl( \frac{b_{ij}}{b_{il}} x \wedge1 \biggr) \,\dif x,\nonumber
\end{eqnarray}
where $\delta_{jk}$ is $1$ if $j = k$ and $0$ if $j \neq k$.
\end{lm}
We illustrate the performance of the M-estimator on two factor models:
a four-dimensional model
with $2$ factors ($p=1\times4=4$), for simulated data sets, and a
three-dimensional
model with $3$ factors ($p=2\times3 = 6$), for real financial data.

The integral on the right-hand side of~(\ref{Elm62}) is to be computed
numerically. For the factor model, the dependence of the matrix
$M(\theta_0)$ on $g$ is too complicated to obtain a general solution
for the optimal function $g$. Since in the previous examples low degree
polynomials gave very good results, and since by the previous lemma
such a choice simplifies the calculations significantly
(numerical integration in dimension $1$, instead of in dimension $d$),
we considered such functions $g$ in a sensitivity analysis. It showed
that the simplest
cases give very good results in terms of root mean squared errors and
that the performance of the M-estimator is quite robust to the
particular choices of $g$. Hence, we suggest using simple, low degree
polynomials for the functions $g$. The functions $g$ in the following
examples are exactly of that type.

\subsection*{Simulation study: Four-dimensional model with two factors}
We simulated $500$ samples of size $n=5000$ from a four-dimensional
model:
\begin{eqnarray*}
X_1 &=& 0.2Z_1 \vee0.8Z_2,
\\
X_2 &=& 0.5Z_1 \vee0.5Z_2,
\\
X_3 &=& 0.7Z_1 \vee0.3Z_2,
\\
X_4 &=& 0.9Z_1 \vee0.1Z_2
\end{eqnarray*}
with independent standard Fr\'echet factors $Z_1$ and
$Z_2$. We have $\theta=(0.2, 0.5,\break 0.7, 0.9)$.

%
\begin{figure}

\includegraphics{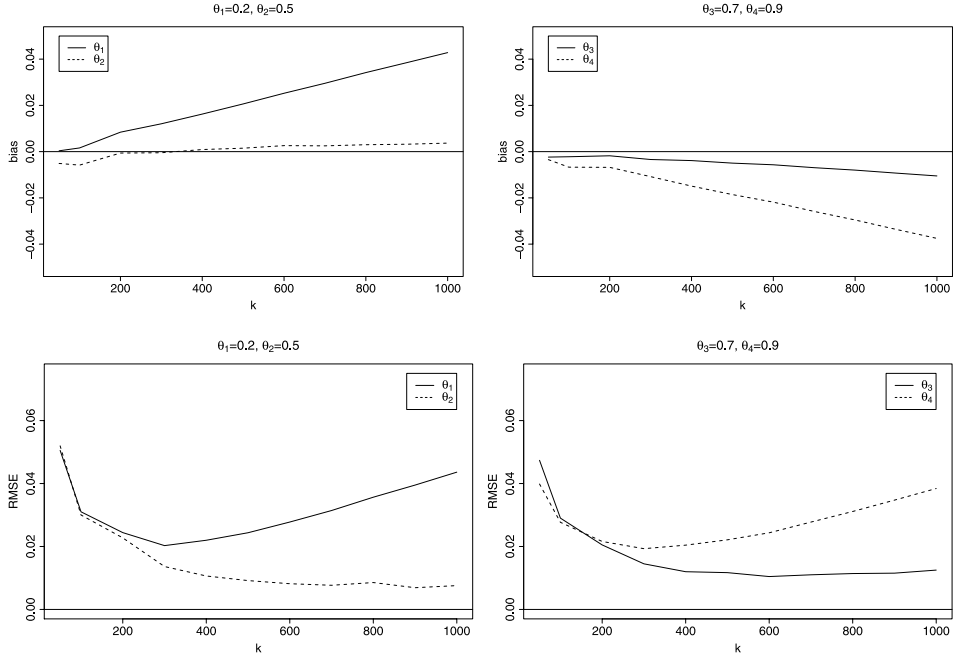}

\caption{Four-dimensional 2-factor model, estimation of $\theta=(0.2,
0.5, 0.7, 0.9)$.}
\label{Figd4r2theta2}\vspace*{-3pt}
\end{figure}

%
%
\begin{figure}[b]
\begin{tabular}{@{}c@{\hspace*{6pt}}c@{}}

\includegraphics{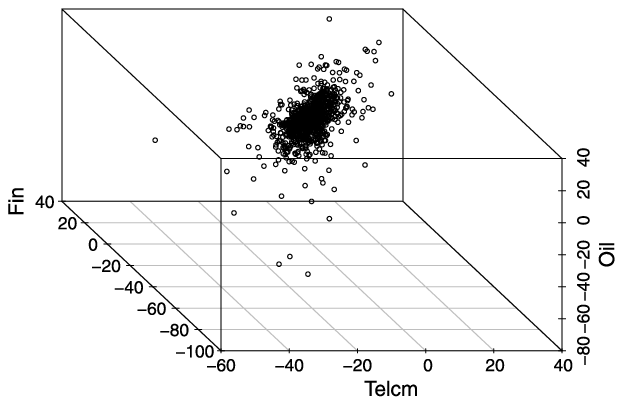}
 & \includegraphics{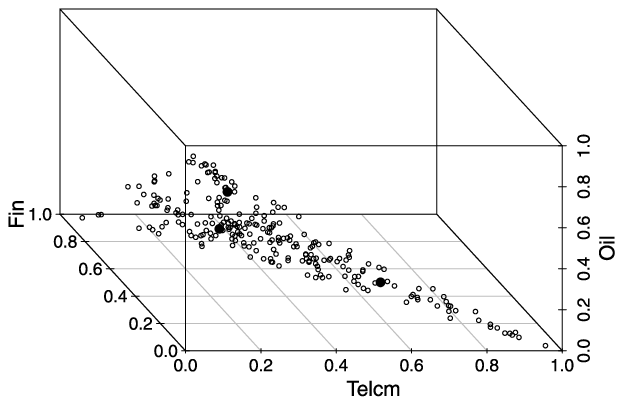}\\
(a) & (b)\vspace*{-3pt}
\end{tabular}
\caption{\textup{(a)} Scatterplot of the original data;
\textup{(b)} plot of the pseudo-data and the three centers.}
\label{Figs3dplots}
\end{figure}

In Figure~\ref{Figd4r2theta2} we show the
bias and the RMSE of the M-estimator based on $q=5$ moment
equations, with auxiliary functions $g_i(x)=x_i$, for $i=1,2,3,4$
and $g_5\equiv1$.
The M-estimator performs very well. For relatively small $k$, the four\vadjust{\goodbreak}
components of $\theta$ are estimated equally well, whereas for larger $k$
the estimator performs somewhat better for parameter values in the
``middle'' of the interval $(0,1)$ than for values near $0$ or $1$.\vspace*{-3pt}

\subsection*{Real data: Three-dimensional model with three factors}
We consider\break monthly negative returns (losses) of three industry
portfolios (Telecommunications, Finance and Oil) over the period July
1, 1926, until December 31, 2009. See Figure~\ref{Figs3dplots}(a) for
the scatterplot of the data; the sample size $n=1002$. The data
are\vadjust{\goodbreak}
available at
\texttt{\href{http://mba.tuck.dartmouth.edu/pages/faculty/ken.french}{http://mba.tuck.dartmouth.edu/pages/}
\href{http://mba.tuck.dartmouth.edu/pages/faculty/ken.french}{faculty/ken.french}}.
We are interested in modeling the losses by a
factor model. In the asset pricing literature [see, e.g., Fama and French
(\citeyear{FF93,FF96})], it is common to model the returns by linear
factor models of type~(\ref{EsumFM}), with three underlying economic
factors. Based on that line of literature, we also consider a
three-factor model for the tails of the three industry portfolios
above; see also \citet{Kleibergen}.

To estimate the parameter vector with $p=2\times3=6$ components, we
need to find a minimum of a 6-dimensional nonlinear criterion function.
To solve such a difficult minimization problem, it is
important to have good starting values. We find a starting parameter
vector by applying the 3-means clustering algorithm [see, e.g.,
\citet{Pollard84}, page 9] to the following pseudo-data:
we transform the data (Telcm, Fin, Oil) to
\[
\bigl(n/(n+1-R_{Ti}), n/(n+1-R_{Fi}), n/(n+1-R_{Oi})
\bigr),\qquad i=1, \ldots, n,
\]
where $R_{Ti}$, $R_{Fi}$ and $R_{Oi}$ are the ranks of the components
of the $i$th observation. Only the entries such that the sum of
their values is greater than the threshold $n/75$ are taken into
account, and subsequently
normalized so that they belong to the unit simplex $\Delta_{3-1}$; see
Figure~\ref{Figs3dplots}(b).
We compute the 3-means cluster centers for these data. Using
equation~(\ref{atoms}), we compute from these three centers
the 6-dimensional starting parameter [as described below
equation~(\ref{E4momentH})] for the minimization routine.
For the criterion function we use $q=7$ functions $g_i$ as follows:
$g_i(x)=x_i$ for $i=1,2,3$, $g_i(x)=x_{i-3}^2$ for $i=4,5,6$,
and $g_7\equiv1$. For different choices of $k$, we obtain the
estimates presented in Table~\ref{TabM-AlogLossALAE2}. For each $k$, we estimate the
loading of the first two factors. This corresponds to the first two
columns of estimated ${b_{ij}}$ for each~$k$. The third columns follow
from the conditions in~(\ref{E4momentH}).

\begin{table}[b]
\caption{Estimates for the factor loadings $b_{ij}$ in
the three-factor model fitted to the tail of the
\mbox{Telcm/Fin/Oil data}}\label{TabM-AlogLossALAE2}
\begin{tabular*}{\tablewidth}{@{\extracolsep{\fill}}lccccc@{}}
\hline
\multicolumn{3}{@{}c}{$\bolds{k=60}$} & \multicolumn{3}{c@{}}{$\bolds{k=90}$} \\
\hline
$0.394$&$0.593$&$0.013$ & $0.344$&$0.616$&$0.040$ \\
$0.691$&$0.211$&$0.098$& $0.701$&$0.216$&$0.083$ \\
$0.358$&$0.062$&$0.580$ & $0.368$&$0.052$&$0.580$\\
\hline
\multicolumn{3}{@{}c}{$\bolds{k=120}$} & \multicolumn{3}{c@{}}{$\bolds{k=150}$} \\
\hline
$0.387$&$0.586$&$0.027$ & $0.388$&$0.581$&$0.031$ \\
$0.695$&$0.215$&$0.090$ & $0.699$&$0.211$&$0.090$ \\
$0.348$&$0.058$&$0.594$ & $0.364$&$0.086$&$0.550$\\
\hline
\end{tabular*}
\end{table}

Observe that the estimates do hardly depend on the choice of $k$. We
see that all three portfolios load substantially on the first factor
(the first column of estimated coefficients, for each $k$),
but Telecommunications loads more on the second\vadjust{\goodbreak} factor (the first lines
of estimated coefficients), and Oil more on the third factor (the third
lines of estimated coefficients).
This indicates that even for only these three portfolios, three factors
are required.

\section{Proofs}\label{sec7}
\label{S7Proofs}

The asymptotic properties of the nonparametric estimator $\hat{l}_n$
are required for the proofs of the asymptotic properties of the
M-estimator~$\hat{\theta}_n$. Consistency of $\hat{l}_n$ [see
(\ref{E4Ln-ln})] for dimension $d=2$ was shown in \citet{Huang92};
cf. \citet{DH98}. In particular, it holds that for every $T>0$,
as $n\to\infty$, $k\to\infty$ and $k/n\to0$,
\[
\sup_{(x_1,x_2)\in[0,T]^2}\bigl|\hat{l}_n(x_1,x_2)-l(x_1,x_2)\bigr|
\pto0.
\]
The proof translates straightforwardly to general dimension $d$, and
together with integrability of $g$ yields consistency of $\tint g
\hat{l}_n$ for $\tint g l = \varphi(\theta_0)$. For the proof of
Theorem~\ref{Tcons}, a technical result is needed.

Let $\mathcal{H}_{k,n}(\theta)\in\R^{p\times p}$ denote the
Hessian matrix of $Q_{k,n}$ as a function of~$\theta$. Let $\mathcal
{H}(\theta)$
be the deterministic, symmetric $p \times p$ matrix whose $(i,j)$th
element, $i, j \in\{1, \ldots, p\}$, is equal to
\begin{eqnarray*}
\bigl(\mathcal{H}(\theta) \bigr)_{ij}&=& 2 \biggl(\frac{\partial}{\partial
\theta_i}
\varphi(\theta) \biggr)^T \biggl(\frac{\partial}{\partial\theta
_j}\varphi(\theta)
\biggr)\\[-2pt]
&&{} -2 \biggl(\frac{\partial^2}{\partial\theta_i\,\partial\theta
_j}\varphi(\theta) \biggr)^T \bigl(
\varphi(\theta_0)-\varphi(\theta) \bigr).
\end{eqnarray*}

%
\begin{lm}\label{Tlemmahessian}
If $k/n\to0$ and if the assumptions of Theorem~\ref{Tcons}\textup{(ii)}
are satisfied, then as $n\to\infty$ and $k\to\infty $, on some closed
neighborhood of $\theta_0$,
\begin{longlist}[(ii)]
\item[(i)] $\mathcal{H}_{k,n}(\theta)\pto\mathcal{H}(\theta)$
uniformly in $\theta$, and
\item[(ii)]$ \P( \mathcal{H}_{k,n}(\theta) \mbox{ is
positive definite} )\to
1$.
\end{longlist}
\end{lm}

\begin{pf}
(i) The Hessian matrix
of $Q_{k,n}$ in $\theta$ is a
$p\times p$ matrix $\mathcal{H}_{k,n}(\theta)$ with elements
$ (\mathcal{H}_{k,n}(\theta) )_{ij} = \partial^2
Q_{k,n}(\theta) / \partial\theta_j \,\partial\theta_i$, for $i,j \in
\{1,\ldots,p\}$, given by
\begin{eqnarray*}
\bigl(\mathcal{H}_{k,n}(\theta) \bigr)_{ij}
&=&2\sum_{m=1}^q\int_{[0,1]^d}g_m(x) \,\frac{\partial}{\partial
\theta_j}l(x;\theta) \,\dif x \cdot\int
_{[0,1]^d}g_m(x)\,\frac{\partial}{\partial\theta_i} l(x;\theta) \,\dif x
\\[-2pt]
&&{}-2\sum_{m=1}^q \int
_{[0,1]^d}g_m(x)\,\frac{\partial^2}{\partial\theta_j\,\partial
\theta_i} l(x;\theta) \,\dif x\\[-2pt]
&&\hspace*{36.2pt}{}\times\int_{[0,1]^d} g_m(x) \bigl(\hat
{l}_n(x)-l(x;\theta)\bigr) \,\dif x
\\[-2pt]
&=& 2 \biggl(\frac{\partial}{\partial\theta_i} \varphi(\theta) \biggr
)^T \biggl(
\frac{\partial}{\partial\theta_j}\varphi(\theta) \biggr) -2 \biggl
(\frac{\partial^2}{\partial\theta_i\,\partial\theta_j}\varphi(
\theta) \biggr)^T
\\[-2pt]
&&{} \times\biggl(\int_{[0,1]^d} g(x) \hat{l}_n(x) \,\dif
x-\varphi(\theta) \biggr).\vadjust{\goodbreak}
\end{eqnarray*}
The consistency of $\tint g\hat{l}_n$ for $\varphi(\theta_0)$
implies
\begin{eqnarray*}
\bigl(\mathcal{H}_{k,n}(\theta) \bigr)_{ij} &\pto& 2 \biggl(
\frac{\partial}{\partial\theta_i}\varphi(\theta) \biggr)^T \biggl(\frac
{\partial}{\partial\theta_j}
\varphi(\theta) \biggr) \\
&&{}-2 \biggl(\frac{\partial^2}{\partial\theta
_i\,\partial\theta_j}\varphi(\theta)
\biggr)^T\bigl(\varphi(\theta_0)-\varphi(\theta)\bigr)
\\
&=& \bigl(\mathcal{H}(\theta) \bigr)_{ij}.
\end{eqnarray*}
Since we assumed that there exists $\eps_0>0$ such that the set
$\{\theta\in\Theta\dvtx\|\theta-\theta_0\|\leq\varepsilon_0\}
=:B_{\eps_0}(\theta_0)$
is closed and thus compact, and since $\varphi$ is assumed to be twice
continuously differentiable, the second
derivatives of $\varphi$ are uniformly bounded on
$B_{\eps_0}(\theta_0)$ and, hence, the convergence above is
uniform on $B_{\eps_0}(\theta_0)$.

(ii) For $\theta=\theta_0$ we get
\[
\bigl(\mathcal{H}(\theta_0) \bigr)_{ij}= 2 \biggl(
\frac{\partial}{\partial\theta_i}\varphi(\theta) \bigg|_{\theta=\theta_0}
\biggr)^T \biggl(
\frac{\partial}{\partial\theta_j}\varphi(\theta) \bigg|_{\theta=\theta_0}
\biggr),
\]
that is,
\[
\mathcal{H}(\theta_0)=2\dot{\varphi}(\theta_0)^T
\dot{\varphi}(\theta_0).
\]
Since $\dot{\varphi}(\theta_0)$ is assumed to be of full rank,
$\mathcal{H}(\theta_0)$ is positive definite. For $\theta$ close to
$\theta_0$, $\mathcal{H}(\theta)$ is also positive definite. Due to
the uniform convergence of $\mathcal{H}_{k,n}(\theta)$ to
$\mathcal{H}(\theta)$ on $B_{\eps_0}(\theta_0)$, the matrix
$\mathcal{H}_{k,n}(\theta)$ is also positive definite on
$B_{\eps_0}(\theta_0)$ with probability tending to one.
\end{pf}

\begin{pf*}{Proof of Theorem \protect\ref{Tcons}}
(i) Fix $\eps> 0$ such that
$0<\varepsilon\leq\varepsilon_0$. Since $\varphi$ is a
homeomorphism, there exists $\delta
> 0$ such that $\theta\in\Theta$ and $\| \varphi(\theta) -
\varphi(\theta_0) \| \leq\delta$ implies $\| \theta- \theta_0 \|
\leq\eps$. In other words, for every $\theta\in\Theta$ such that
$\| \theta- \theta_0 \| > \eps$, we have $\| \varphi(\theta) -
\varphi(\theta_0) \| > \delta$. Hence, on the event
\[
A_n = \biggl\{\biggl\| \varphi(\theta_0) - \tint g
\hat{l}_n \biggr\| \leq\delta/2\biggr\}
\]
for every $\theta\in\Theta$ with $\| \theta- \theta_0 \| > \eps$,
necessarily,
\begin{eqnarray*}
\biggl\| \varphi(\theta) - \tint g \hat{l}_n \biggr\|
&\geq&\bigl\| \varphi(\theta) -
\varphi(\theta_0) \bigr\| - \biggl\| \varphi(\theta_0) - \tint g
\hat{l}_n \biggr\|\\
&>& \delta- \delta/2 = \delta/2 \geq\biggl\| \varphi(
\theta_0) - \tint g \hat{l}_n \biggr\|.
\end{eqnarray*}
As a consequence, on the event $A_n$, we have
\[
\inf_{\theta\dvtx\| \theta- \theta_0 \| > \eps} \biggl\| \varphi(\theta) -
\tint
g \hat{l}_n \biggr\| >
\min_{\theta\dvtx\| \theta- \theta_0 \| \leq\eps} \biggl\| \varphi(\theta) -
\tint g \hat{l}_n \biggr\|,
\]
where we can write the minimum on the right-hand side since the set $\{
\theta\in\Theta\dvtx\| \theta- \theta_0 \| \leq\eps\}$ is
closed and thus compact for $0\leq\eps\leq\eps_0$. Hence, on the
event~$A_n$, the ``argmin'' set $\hat{\Theta}_n$ is
nonempty and is contained in the closed ball of radius $\eps$
centered at $\theta_0$. Finally, $\P(A_n) \to1$ by weak consistency
of $\int g \hat{l}_n$ for $\int g l = \varphi(\theta_0)$.

(ii) In the proof of (i) we have seen that, with probability tending
to one, the proposed M-estimator exists and it is contained in a
closed ball around $\theta_0$. In Lem\-ma~\ref{Tlemmahessian} we
have shown that the criterion function is, with probability tending
to one, strictly convex on such a closed ball around $\theta_0$ and,
hence, with probability tending to one, the minimizer of the
criterion function is unique.
\end{pf*}

For $i=1,\ldots,n$ let
\[
U_i:= (U_{i1},\ldots,U_{id}):=
\bigl(1-F_1(X_{i1}),\ldots,1-F_d(X_{id})
\bigr)
\]
and denote
\begin{eqnarray*}
Q_{nj}(u_j)&:=&U_{\lceil nu_j \rceil\dvtx n,j},\qquad j=1,\ldots,d,
\\
S_{nj}(x_j)&:=&\frac{n}{k}Q_{nj} \biggl(
\frac{kx_j}{n} \biggr),\qquad j=1,\ldots,d,
\\
S_n(x)&:=&\bigl(S_{n1}(x_1),
\ldots,S_{nd}(x_d)\bigr),
\end{eqnarray*}
where $U_{1\dvtx n,j}\!\leq\!\cdots\!\leq\! U_{n\dvtx n,j}$ are the order
statistics of
$U_{1j},\ldots,U_{nj}$, $j\!=\!1,\ldots,d$, and $\lceil a \rceil$ is the
smallest integer not smaller than $a$. Write
\begin{eqnarray*}
V_n(x) :\!&=& \frac{n}{k}\P\biggl( U_{11}\leq
\frac{kx_1}{n} \mbox{ or }\ldots\mbox{ or } U_{1d}\leq
\frac{kx_d}{n} \biggr),
\\
T_n(x) :\!&=& \frac{1}{k}\sum_{i=1}^n
\one\biggl\{ U_{i1}< \frac
{kx_1}{n} \mbox{ or }\ldots\mbox{ or }
U_{id}< \frac{kx_d}{n} \biggr\},
\\
\hat{L}_n(x) :\!&=& \frac{1}{k}\sum
_{i=1}^n \one\biggl\{ U_{i1}<
\frac
{k}{n}S_{n1}(x_1) \mbox{ or }\ldots\mbox{ or }
U_{id}< \frac
{k}{n}S_{nd}(x_d) \biggr\}
\\
& =& \frac{1}{k}\sum_{i=1}^n\one
\bigl\{R_{i}^1>n+1-kx_1 \mbox{ or }\ldots\mbox{
or } R_{i}^d>n+1-kx_d \bigr\}
\end{eqnarray*}
and note that
\[
\hat{L}_n(x) = T_n\bigl(S_n(x)\bigr).
\]
With probability one, for every $x$ and for every $j\in\{1, \ldots,
d\}$, there is at most one $i$ such
that $n+\frac{1}{2}-kx_j<R_i^j\leq n+1-kx_j$. Hence,
%
%
\begin{equation}
\label{E4Ln-ln} \sup_{x\in[0,1]^d}\sqrt{k}\bigl\llvert\hat{l}_n(x)-
\hat{L}_n(x)\bigr\rrvert\leq\frac{d}{\sqrt{k}}\to0.
\end{equation}
This shows that the asymptotic properties of $\hat{l}_n$ and $\hat
{L}_n$ are the
same. With the notation $v_n(x)=\sqrt{k}(T_n(x)-V_n(x))$, we have
the following result.

%
\begin{pr} \label{Tprop}
Let $T>0$ and denote $A_x:=\{u\in[0,\infty]^d\dvtx u_1\leq x_1$
or $\ldots$ or $u_d\leq x_d\}$. There exists a sequence of
processes $\tilde{v}_n$ such that, for all~$n$, $\tilde{v}_n\eqd v_n$
and there exists a Wiener process $W_l(x):=W_\Lambda(A_x)$ such that
as $n\to\infty$,
%
%
\begin{equation}
\label{E4proposition} \sup_{x\in[0,2T]^d}\bigl|\tilde{v}_n(x)-W_l(x)\bigr|
\pto0.
\end{equation}
\end{pr}

The result follows from Theorem 3.1 in \citet{E97}. From the proofs
there it follows that a single Wiener process, instead of the
sequence in the original statement of the theorem, can be used, and
that convergence holds almost surely, instead of in probability,
once the Skorohod construction is introduced. From now on, we work
on this new (Skorohod) probability space, but keep the old notation,
without the tildes. In particular, we have convergence of the
marginal processes:
\[
\sup_{x_j\in[0,2T]}\bigl|v_{nj}(x)-W_{l,j}(x_j)\bigr|\to0
\qquad\mbox{a.s.},\qquad j=1,\ldots,d,
\]
where $v_{nj}(x_j):=v_n((0,\ldots,0,x_j,0,\ldots,0))$. The
\citet{Vervaat72} lemma implies
%
%
\begin{equation}
\label{E4Snx-x} \sup_{x_j\in[0,2T]}\bigl|\sqrt{k}\bigl(S_{nj}(x_j)-x_j
\bigr)+W_{l,j}(x_j)\bigr|\to0 \qquad\mbox{a.s.},\qquad
j=1,\ldots,d.\hspace*{-28pt}
\end{equation}

\begin{pf*}{Proof of Theorem~\ref{TanL}}
Write
\begin{eqnarray*}
&&
\sqrt{k} \bigl(\hat{L}_n(x)-l(x) \bigr)
\\
&&\qquad= \sqrt{k} \bigl( T_n\bigl(S_n(x)\bigr)-V_n
\bigl(S_n(x)\bigr) \bigr) + \sqrt{k} \bigl( V_n
\bigl(S_n(x)\bigr)-l\bigl(S_n(x)\bigr) \bigr)\\
&&\qquad\quad{} + \sqrt{k}
\bigl( l\bigl(S_n(x)\bigr)-l(x) \bigr)
\\
&&\qquad=: D_1(x) + D_2(x) + D_3(x).
\end{eqnarray*}

\textsc{Proof of} $\sup_{x\in[0,T]^d}|D_1(x)-W_l(x)|\pto0$.\quad
We have
\[
D_1(x)=\sqrt{k} \bigl(T_n\bigl(S_n(x)
\bigr)-V_n\bigl(S_n(x)\bigr) \bigr)=v_n
\bigl(S_n(x)\bigr).
\]
It holds that
\begin{eqnarray*}
&& \sup_{x\in[0,T]^d}\bigl|D_1(x)-W_l(x)\bigr|
\\
&&\qquad \leq \sup_{x\in[0,T]^d}\bigl\llvert D_1(x) - W_l
\bigl(S_n(x)\bigr)\bigr\rrvert\\
&&\qquad\quad{}+ \sup_{x\in[0,T]^d}\bigl\llvert
W_l\bigl(S_n(x)\bigr)-W_l(x)\bigr\rrvert.
\end{eqnarray*}
Because of~(\ref{E4Snx-x}), this is, with probability tending to
one, less than or equal to
\[
\sup_{y\in[0,2T]^d}\bigl\llvert v_n(y)-W_l(y)\bigr
\rrvert+ \sup_{x\in
[0,T]^d}\bigl\llvert W_l\bigl(S_n(x)
\bigr)-W_l(x)\bigr\rrvert.
\]
Both terms tend to zero in probability, the first one by Proposition
\ref{Tprop}, the second one because of the uniform continuity of
$W_l$ and~(\ref{E4Snx-x}).

\textsc{Proof of} $\sup_{x\in[0,T]^d}|D_2(x)|\pto0$.\ \
Because of~(\ref{E4Snx-x}), with probability~tending to one,
$\sup_{x\in[0,T]^d}|D_2(x)|$ is less\vspace*{1pt} than or equal to
$\sup_{y\in[0,2T]^d}\sqrt{k}|V_n(y)-l(y)|$, which in turn, because
of conditions (C1) and (C2), is equal to
\[
\sqrt{k}O \biggl( \biggl(\frac{k}{n} \biggr)^\alpha\biggr) = O
\biggl( \biggl(\frac{k}{n^{2\alpha/(1+2\alpha)}} \biggr)^{
1/2+\alpha} \biggr) = o(1).
\]

\textsc{Proof of}\vspace*{1pt} $\sup_{x\in[0,T]^d}|D_3(x)+\sum
_{j=1}^dl_j(x)W_{l,j}(x_j)|\pto0$.\quad
Due to the existence of the first derivatives, we can use the mean
value theorem to write
\[
\frac{1}{\sqrt{k}}D_3(x) = l \bigl(S_n(x) \bigr) - l(x)
=\sum_{j=1}^d\bigl(S_{nj}(x_j)-x_j
\bigr)\cdot l_j(\xi_n)
\]
with $\xi_n$ between $x$ and $S_{n}(x)$. Therefore,
\begin{eqnarray*}
&&
\sup_{x\in[0,T]^d}\Biggl|D_3(x)+\sum_{j=1}^dl_j(x)W_{l,j}(x_j)\Biggr|\\
&&\qquad\leq\sum_{j=1}^d\bigl|l_j(
\xi_n)\sqrt{k}\bigl(S_{nj}(x_j)-x_j
\bigr)+l_j(x)W_{l,j}(x_j)\bigr|.
\end{eqnarray*}
Note that all the terms on the right-hand side of the above
inequality can be dealt with in the same way. Therefore, we consider
only the first
term. For $\delta\in(0,T)$, this term is bounded by
\begin{eqnarray*}
&&\sup_{x\in[0,T]^d} \bigl|l_1(\xi_n)\bigr| \cdot
\sup_{x_1\in[0,T]}\bigl| \sqrt{k}\bigl(S_{n1}(x_1)-x_1
\bigr)+W_{(l,1)}(x_1)\bigr|
\\
&&\quad{} + \sup_{x\in[\delta, T]\times[0,T]^{d-1}}\bigl|l_1(\xi_n)-l_1(x)\bigr|
\cdot\sup_{x_1\in[0,T]}\bigl|W_{(l,1)}(x_1)\bigr|
\\
&&\quad{} + \sup_{x\in[0, \delta]\times[0,T]^{d-1}}\bigl|l_1(\xi_n)-l_1(x)\bigr|
\cdot\sup_{x_1\in[0,\delta]}\bigl|W_{(l,1)}(x_1)\bigr|
\\
&&\qquad =: D_4\cdot D_5 + D_6\cdot
D_7+ D_8\cdot D_9.
\end{eqnarray*}
Observe that $0\leq l_1\leq1$. Also, since $l_1$ is continuous on
$[\delta/2, T]\times[0,T]^{d-1}$, it is uniformly continuous on that
region. We have $D_5\pto0$ by~(\ref{E4Snx-x}), so $D_4\cdot
D_5\pto0$. The uniform continuity
of $l_1$ and the fact that almost surely $D_7<\infty$ yield
$D_6\cdot D_7\pto0$. Finally, for every $\varepsilon>0$, we can find
a $\delta$ such that, with probability at least $1-\varepsilon$,
$D_9<\varepsilon$
and, hence, $D_8\cdot D_9<\varepsilon$.

Applying~(\ref{E4Ln-ln}) completes the proof. 
\end{pf*}

%
\begin{pr} \label{TANintegrals}
If conditions \textup{(C1)} and \textup{(C2)} from Theorem
\ref{Tantheta} hold, then as $n\to\infty$,
%
%
\begin{equation}
\label{E4ANintegrals} \sqrt{k}\int_{[0,1]^d}g(x) \bigl(
\hat{l}_n(x)-l(x) \bigr)\,\dif x\dto\tilde{B}.
\end{equation}
\end{pr}

\begin{pf}
Throughout the proof we
write $l(x)$ instead of
$l(x;\theta_0)$. Also, since $l$ does not need to be differentiable,
we will use notation $l_j(x)$, $j=1,\ldots,d$, to denote the
right-hand partial derivatives here. Let $D_1(x), D_2(x),\break D_3(x)$
be as in the proof of Theorem~\ref{TanL} and take $T=1$. Then
\begin{eqnarray*}
\hspace*{-4pt}&&\biggl|\sqrt{k} \biggl(\int_{[0,1]^d}g(x)\hat{L}_n(x)\,\dif x
- \int_{[0,1]^d}g(x)l(x)\, \dif x \biggr)-\tilde{B} \biggr|
\\
\hspace*{-4pt}&&\qquad\leq\sup_{x\in[0,1]^d}\bigl|D_1(x)-W_l(x)\bigr|\int
_{[0,1]^d}\bigl|g(x)\bigr|\,\dif x
+ \sup_{x\in[0,1]^d}\bigl|D_2(x)\bigr|
\int_{[0,1]^d}\bigl|g(x)\bigr|\,\dif x
\\
\hspace*{-4pt}&&\qquad\quad{}+ \int_{[0,1]^d}\bigl|g(x,y)\bigr|\cdot\Biggl\llvert D_3(x) +
\sum_{j=1}^dl_j(x)W_{l,j}(x_j)
\Biggr\rrvert\,\dif x.
\end{eqnarray*}
The first two terms on the right-hand side converge to zero in
probability due to integrability of $g$ and uniform convergence of
$D_1(x)$ and $D_2(x)$, which was shown in the proof of Theorem
\ref{TanL}. The third term needs to be treated separately, as the
condition on continuity (and existence) of partial derivatives is no
longer assumed to hold.

Let $\omega$ be a point in the Skorohod probability space introduced
before the proof of Theorem~\ref{TanL} such that for all
$j=1,\ldots,d$,
\[
\sup_{x_j\in[0,1]}\bigl|W_{l,j}(x_j)\bigr|<+\infty
\quad\mbox{and}\quad
\sup_{x_j\in[0,1]}\bigl|\sqrt{k}\bigl(S_{nj}(x_j)-x_j
\bigr)+W_{l,j}(x_j)\bigr|\to0.
\]
For such $\omega$ we will show by means of dominated convergence
that
%
%
\begin{equation}
\label{E43rdterm}
\quad\int_{[0,1]^d}\bigl|g(x)\bigr|\cdot\Biggl\llvert\sqrt{k}
\bigl(l\bigl(S_n(x)\bigr)-l(x) \bigr) +\sum
_{j=1}^dl_j(x)W_{l,j}(x_j)
\Biggr\rrvert\,\dif x \to0.
\end{equation}

\textsc{Proof of the pointwise convergence.}\quad If $l$ is differentiable,
convergence of the above integrand to zero follows from the
definition of partial derivatives and~(\ref{E4Snx-x}). Since
this might fail only on a set of Lebesgue measure zero, the
convergence of the integrand to zero holds almost everywhere on
$[0,1]^d$.\vadjust{\goodbreak}

\textsc{Proof of the domination.}\quad Note that from expressions for
(one-sided) partial derivatives~(\ref{E4lgrad}), and the moment
conditions~(\ref{E4moment}), it follows that $0 \leq l_j(x) \leq1$, for
all $x\in[0,1]^d$ and all $j=1,\ldots,d$.

We get
\begin{eqnarray*}
&&
\bigl|g(x)\bigr| \cdot\Biggl\llvert\sqrt{k} \bigl(l \bigl(S_n(x) \bigr)-l(x)
\bigr)+\sum_{j=1}^dl_j(x)W_{l,j}(x_j)
\Biggr\rrvert
\\
&&\qquad\leq\bigl|g(x)\bigr|\cdot\Biggl(\sqrt{k}\bigl|l\bigl(S_n(x)\bigr)-l(x)\bigr|+\sum
_{j=1}^d\bigl|W_{l,j}(x_j)\bigr|
\Biggr).
\end{eqnarray*}
Using the definition of function $l$ and uniformity of
$1-F_j(X_{1j})$, we have, for all $j=1,\ldots,d$,
\[
\bigl|l(x_1, \ldots,x_{j-1}, x_j, x_{j+1},
\ldots,x_d)-l\bigl(x_1, \ldots,x_{j-1},
x_j', x_{j+1},\ldots,x_d\bigr)\bigr|
\leq\bigl|x_j-x_j'\bigr|.
\]
Hence, we can write
\begin{eqnarray*}
&&
\sup_{x\in[0,1]^d}\sqrt{k}\bigl|l\bigl(S_n(x)\bigr)-l(x)\bigr| \\
&&\qquad \leq
\sup_{x\in[0,1]^d}\sqrt{k}\bigl|l\bigl(S_n(x)\bigr)-l
\bigl(x_1, S_{n2}(x_2),\ldots,S_{nd}(x_d)
\bigr)\bigr|
\\
&&\qquad\quad{} + \sup_{x\in[0,1]^d}\sqrt{k}\bigl|l\bigl(x_1, S_{n2}(x_2),S_{n3}(x_3),
\ldots,S_{nd}(x_d)\bigr)
\\
&&\qquad\quad\hspace*{74pt}{} - l\bigl(x_1,x_2, S_{n3}(x_3),
\ldots,S_{nd}(x_d)\bigr)\bigr|
+ \cdots
\\
&&\qquad\quad{} + \sup_{x\in[0,1]^d}\sqrt{k}\bigl|l\bigl(x_1,\ldots,x_{d-1},
S_{nd}(x_d)\bigr)-l(x)\bigr|
\\
&&\qquad \leq \sum_{j=1}^d \sup_{x_j\in[0,1]}
\sqrt{k}\bigl|S_{nj}(x_j)-x_j\bigr| = O(1).
\end{eqnarray*}
Since for all\vspace*{2pt} $j=1,\ldots,d$ we have
${\sup_{x_j\in[0,1]}}|W_{l,j}(x_j)|<+\infty$, the proof of (\ref
{E43rdterm}) is complete. This, together with~(\ref{E4Ln-ln}), finishes
the proof of the proposition.
\end{pf}

Let $\nabla Q_{k,n}(\theta)\in\R^{p\times1}$ be the gradient vector
of $Q_{k,n}$ at $\theta$. Put
\[
V(\theta):= 4 \dot{\varphi}(\theta)^T \Sigma(\theta) \dot{\varphi}(
\theta)\in\R^{p\times p}.
\]

%
\begin{lm}\label{Tlemmagradient}
If the assumptions of Theorem~\ref{Tantheta} are satisfied, then as
$n\to\infty$,
\[
\sqrt{k}\nabla Q_{k,n}(\theta_0)\dto N\bigl(0,V(
\theta_0)\bigr).
\]
\end{lm}

\begin{pf}
The gradient vector of
$Q_{k,n}$ with respect to
$\theta$ in $\theta_0$ is
\[
\nabla Q_{k,n}(\theta_0) = \biggl(\frac{\partial}{\partial\theta_1}
Q_{k,n}(\theta) \bigg|_{\theta=\theta_0},\ldots,\frac{\partial
}{\partial\theta_p}
Q_{k,n}(\theta) \bigg|_{\theta=\theta_0} \biggr)^T,
\]
where for $i=1,\ldots,p$,
\begin{eqnarray*}
&&\frac{\partial}{\partial\theta_i} Q_{k,n}(\theta) \bigg|_{\theta=\theta
_0}\\
&&\qquad=-2\sum
_{m=1}^q\int_{[0,1]^d}g_m(x)\,
\frac{\partial}{\partial\theta_i} l(x;\theta) \bigg|_{\theta=\theta_0}\,\dif
x\\
&&\hspace*{0pt}\qquad\quad{}\times\int
_{[0,1]^d}g_m(x) \bigl(\hat{l}_n(x)-l(x;
\theta_0)\bigr)\,\dif x.
\end{eqnarray*}
Using vector notation, we obtain
\[
\nabla Q_{k,n}(\theta_0) = -2\dot{\varphi}(
\theta_0)^T\cdot\int_{[0,1]^d} g(x)
\bigl(\hat{l}_n(x) - l(x;\theta_0)\bigr) \,\dif x.
\]
Equation~(\ref{E4Ln-ln}) and the proof of Proposition~\ref{TANintegrals}
imply that
\begin{eqnarray*}
\sqrt{k}\nabla Q_{k,n}(\theta_0) = -2\dot{\varphi}(
\theta_0)^T\cdot\int_{[0,1]^d} g(x)
\sqrt{k} \bigl(\hat{l}_n(x) - l(x;\theta_0) \bigr) \,\dif x
\dto-2\dot{\varphi}(\theta_0)^T \tilde{B}.
\end{eqnarray*}
The limit\vspace*{1pt} distribution of $\sqrt{k}\nabla Q_{k,n}(\theta_0)$ is
therefore zero-mean Gaussian with covariance matrix $V(\theta_0) =
4\dot{\varphi}(\theta_0)^T \Sigma(\theta_0)\dot{\varphi}(\theta_0)$.
\end{pf}

\begin{pf*}{Proof of Theorem~\ref{Tantheta}}
Consider the function
$f(t):=\nabla Q_{k,n}(\theta_0+t(\hat{\theta}_n-\theta_0))$,
$t\in[0,1]$. The mean value theorem yields
\[
\nabla Q_{k,n}(\hat{\theta}_n)=\nabla Q_{k,n}(
\theta_0) + \mathcal{H}_{k,n}(\tilde{\theta}_n)
(\hat{\theta}_n-\theta_0)
\]
for some $\tilde{\theta}_n$ between $\theta_0$ and $\hat{\theta}_n$.
First note that, with probability tending to one, $0=\nabla
Q_{k,n}(\hat{\theta}_n)$, which follows from the fact that
$\hat{\theta}_n$ is a minimizer of $Q_{k,n}$ and that, with
probability tending to one, $\hat{\theta}_n$ is in an open ball
around $\theta_0$. By the consistency of $\hat{\theta}_n$, we have
that $\tilde{\theta}_n\pto\theta_0$, and since the convergence of
$\mathcal{H}_{k,n}$ to $\mathcal{H}$ is uniform on a neighborhood of
$\theta_0$, we get that $\mathcal{H}_{k,n}(\tilde{\theta}_n)\pto
\mathcal{H}(\theta_0)$. Hence, $\sqrt{k}(\hat{\theta}_n-\theta_0)\dto
N(0,M(\theta_0))$.
\end{pf*}

\begin{pf*}{Proof of Corollary~\ref{Tcorollary}}
As in Lemma 7.2 in
\citet{EKS08}, we can see that if $\theta\mapsto H_\theta$ is weakly
continuous at $\theta_0$, then $\theta\mapsto\Sigma(\theta)$ is
continuous at $\theta_0$. This, together with the assumption that
$\varphi$ is twice continuously differentiable and $\dot{\varphi
}(\theta_0)$ is of
full rank, yields that $\theta\mapsto V(\theta)$ is continuous at
$\theta_0$.
The above assumption also implies that $\theta\mapsto
\mathcal{H}(\theta)$ is continuous at $\theta_0$, which, with the
positive definiteness of $\mathcal{H}(\theta)$ in a~neighborhood of
$\theta_0$, shows that if $\theta\mapsto H_\theta$ is weakly
continuous at $\theta_0$, then $\theta\mapsto M(\theta)=
\mathcal{H}(\theta)^{-1}V(\theta)\mathcal{H}(\theta)^{-1}$ is
continuous at $\theta_0$. Hence, we obtain
\[
M(\hat{\theta}_n)^{-1/2}\sqrt{k}(\hat{\theta}_n-
\theta_0)\dto N(0,I_p),
\]
which yields~(\ref{Tcorollary}).
\end{pf*}

\begin{pf*}{Proof of Theorem~\ref{TTest}}
Theorem~\ref{Tantheta} and
the arguments used in the proof of Corollary~\ref{Tcorollary} imply
that, as $n\to\infty$,
%
%
\begin{equation}
\label{E4testeq1} M_2^{-1/2}\bigl(\hat{\theta}_1,
\theta_2^*\bigr)\sqrt{k}\bigl(\hat{\theta}_2-
\theta_2^*\bigr)\dto N(0, I_r)
\end{equation}
and hence~(\ref{E4test}).
\end{pf*}

\begin{pf*}{Proof of Lemma~\ref{TlFMl}}
We have
\begin{eqnarray*}
l(x_1, \ldots, x_d) & = & \lim_{t \to\infty} t \P
\bigl(1-F_1(X_1) \leq x_1/t \mbox{ or } \ldots
\mbox{ or } 1-F_d(X_d) \leq x_d/t \bigr)
\\[-1pt]
& = & \lim_{t \to\infty} t \P\bigl(1-F_{Y_1}(Y_1) \leq
x_1/t \mbox{ or } \ldots\mbox{ or } 1-F_{Y_d}(Y_d)
\leq x_d/t \bigr)
\\[-1pt]
& = & \lim_{t \to\infty} t \P\biggl(Y_1 \geq\frac{t\sum_{i=1}^r
a_{i1}^\nu}{x_1}
\mbox{ or } \ldots\mbox{ or } Y_d \geq\frac
{t\sum_{i=1}^r a_{id}^\nu}{x_d} \biggr)
\\[-1pt]
& = & \lim_{t \to\infty} t \P\biggl( \bigcup_{1 \leq j \leq d}
\bigcup_{1 \leq i \leq r} \biggl\{ W_i \geq
\frac{t\sum_{i=1}^r a_{ij}^\nu}{a_{ij}^\nu x_j} \biggr\} \biggr)
\\[-1pt]
& = & \lim_{t \to\infty} t \P\biggl( \bigcup_{1 \leq i \leq r}
\biggl\{ W_i \geq\min_{1
\leq j
\leq d} \frac{t\sum_{i=1}^r a_{ij}^\nu}{a_{ij}^\nu x_j} \biggr\}
\biggr)
\\[-1pt]
& = & \lim_{t \to\infty} t \sum_{i=1}^r
\P\biggl(W_i \geq\min_{1 \leq j \leq d} \frac{t\sum_{i=1}^r a_{ij}^\nu
}{a_{ij}^\nu x_j} \biggr)
\\[-1pt]
& = & \lim_{t \to\infty} \sum_{i=1}^r t
\biggl( 1 - \exp\biggl\{ - \frac{1}{t}\max_{1 \leq j
\leq d}
\frac{a_{ij}^\nu x_j}{\sum_{i=1}^r a_{ij}^\nu} \biggr\} \biggr)
\\[-1pt]
& = & \sum_{i=1}^r \max_{1 \leq j \leq d}
\biggl\{\frac{a_{ij}^\nu
x_j}{\sum_{i=1}^r a_{ij}^\nu} \biggr\} =: \sum_{i=1}^r
\max_{1 \leq j \leq d} \{b_{ij}x_j\}
\end{eqnarray*}
as required.
\end{pf*}

\begin{pf*}{Proof of Lemma~\ref{TlFMintegrals}}
Fix $i \in\{1, \ldots, r\}$. We have
\[
\int_{[0, 1]^d} x_k^s
\max_{1 \leq j \leq d} \{b_{ij} x_j\} \,\dif x = \sum
_{j=1}^d \int_{[0, 1]^d}
x_k^s (b_{ij} x_j) \one\Bigl(
b_{ij} x_j \geq\max_{l
\neq j} \{b_{il}
x_l\} \Bigr) \,\dif x.
\]
Write the integral as a double integral, the outer integral with
respect to $x_j \in[0, 1]$ and the inner integral\vadjust{\goodbreak} with respect to
$x_{-j} = (x_l)_{l \neq j} \in\mathbb{R}^{d-1}$ over the relevant
domain. We find
\[
\int_{[0, 1]^d} x_k^s
\max_{1 \leq j \leq d} \{b_{ij} x_j\} \,\dif x = \sum
_{j=1}^d \int_0^1
b_{ij} x_j \int_{0 < x_l < ({b_{ij}}/{b_{il}}) x_j \wedge
1}
x_k^s \,\dif x_{-j} \,\dif x_j.
\]
After some long, but elementary computations, this simplifies to the
stated expression.
\end{pf*}

\section*{Acknowledgments}

We are grateful to Axel B\"ucher for pointing out that the original
condition (C3) of Theorem~\ref{TanL} was too restrictive. We also
like to thank the Associate Editor and two referees for a thorough
reading of the manuscript and for many thoughtful comments that led to
this improved version.


%

\printaddresses

\end{document}